\documentclass{article}
\usepackage{amsfonts}
\usepackage{graphicx}
\usepackage{xcolor}
\usepackage{mathtools}
\usepackage{hyperref}
\usepackage{amsthm}
\usepackage{physics}
\usepackage{tikz}
\usetikzlibrary{arrows.meta,arrows}
\usepackage{cleveref}
\usepackage{subcaption}

\newtheorem{theorem}{Theorem}
\newtheorem{lemma}{Lemma}

\newtheorem{corollary}{Corollary}
\newtheorem{definition}{Definition}
\newtheorem{remark}{Remark}

\newcommand{\lstar}{\mathrm{star}^{-}}

\newcommand{\R}{\mathbb{R}}

\newcommand{\relu}{\textrm{ReLU}}
\newcommand{\CF}{\mathcal{C}(F)}
\newcommand{\st}{\mathrm{star}}
\newcommand{\lk}{\mathrm{link}}
\newcommand{\llk}{\mathrm{link}^-}
\newcommand{\rk}{\mathrm{rk}}

\DeclareMathOperator{\sgn}{sign}

\title{Combinatorial Regularity for Relatively Perfect Discrete Morse Gradient Vector Fields of ReLU Neural Networks
	}
\author{Robyn Brooks, Marissa Masden}

\begin{document}
	\maketitle
	\tableofcontents

	\newpage 
	
	\section{Introduction}
	
	Much of the recent progress of machine learning is due to exponential growth in available computational power. 
	This growth has enabled neural network models with trillions of parameters to be trained on extremely large datasets, with exceptional results. 
	In contrast, environmental and economic concerns lead to the question of the minimal network architecture necessary to perform a specific task. 
	
	The theoretical characterization of the exact capabilities of fixed neural network architectures is still ongoing.
	To understand the classification capacity of ReLU neural network functions geometrically, researchers have investigated {attributes such as} the expected number of linear regions \cite{HaninRolnick} and the average curvature of decision boundary \cite{curvature}. 
	To describe the topological capacity of such networks, one measure of interest is the achievable Betti numbers of decision regions and decision boundaries for a fixed neural network architecture\cite{ergen2024topological, guss}.  
	
	In this vein, we continue the development of algorithmic tools for characterizing the topological behavior of fully-connected, feedforward ReLU neural network functions using piecewise linear (PL) Morse and discrete Morse theory. 
	Previous work has established that most ReLU neural networks are PL Morse, and identified conditions for a point in the input space of a ReLU neural network to be regular or have PL Morse index $k$ \cite{GLMas}.
	However, discrete Morse theory provides graph-theoretic algorithms for speeding up computations and narrowing down which computations even need to be performed \cite{HERSH2005}. Furthermore, discrete Morse theory has tools for canceling critical cells which allow for the simplification of sets of critical cells.  
	
	The usefulness of discrete Morse theory in algorithmic computations and understanding the topology of cellular spaces motivates us to bridge the existing gap between PL Morse functions and discrete Morse functions in the context of ReLU neural networks.

		\subsection{Contributions and Related Work}
	
	We follow a framework for the polyhedral decomposition of the input space of a ReLU neural network $F$  by the \emph{canonical polyhedral complex} which we will denote as $\CF$, as introduced in \cite{GrigsbyLindsey} and defined in \Cref{def:canonicalpoly}.  
	This framework is expanded in \cite{GLMas}, introducing a PL Morse characterization, and in \cite{Masden}, introducing the combinatorial characterization of the polyhedral decomposition by sign sequences.  
	 
	Unfortunately, computing the canonical polyhedral complex and topological properties of its level and sublevel sets is expensive and can only be done for small neural networks, in part because $\CF$ has exponentially many vertices in the input dimension of $F$. It is this drawback that is the main motivation for the construction of a discrete function associated to $F$; such a function will retain the topological information desired, but should ease issues of hight computational complexity.
	
	We build on existing literature by providing a translation from the PL Morse characterization from \cite{GLMas} to a discrete Morse gradient vector field by exploiting both the combinatorics from \cite{Masden} and a technique for generating relatively perfect discrete gradient vector fields from \cite{Fugacci}.      
    The main result of this paper (\Cref{thm:relativelyperfectDGVF}) provides a canonical way to directly construct a discrete Morse function (defined on the input space) which captures the same topological information as the existing PL Morse function.
	This paper is intended as an intermediary which should allow further computational tools to be developed.
	
	While not the primary goal of this paper, a secondary contribution which we hope to highlight is the development of additional tools for translating from piecewise linear functions on non-simplicial complexes to discrete Morse gradient vector fields. 
	To our knowledge, all current theory relating PL Morse functions to discrete Morse functions on the same complexes is detailed in \cite{Fugacci}.  We also provide realizability results for certain shallow networks (\Cref{thm:all_realizable}) - namely, we are able to characterize the possible homotopy types of the descision boundary for a $(n,n+1,1)$ generic PL Morse ReLU neural network $F$ by restricting the number and possible indices of the critical points of $F$.

	\subsection{Outline}

	In \Cref{sec:polyandmorse}, we review the mathematical tools we are using. In \Cref{sec:canonicalpoly}, we review the construction of the canonical polyhedral complex and share some novel realizability results for shallow neural networks. 
	\Cref{sec:relperfect} contains our main result (\Cref{thm:relativelyperfectDGVF}). This theorem provides a translation from a PL Morse ReLU neural network $F$ to a relatively perfect discrete gradient vector field on $\CF$. 
	In \Cref{sec:computation} we discuss some of the issues surrounding effective computational implementation, and in \Cref{sec:conclusion}, we provide concluding remarks.


	
	\section{Background: Polyhedral Geometry and Morse Theories}
	\label{sec:polyandmorse}

	A ReLU neural network is a piecewise linear function on a polyhedral complex which we call the \textit{canonical polyhedral complex}.  
	Before developing constructions on this specific object, we review some general constructions in polyhedral and piecewise linear geometry which we will use repeatedly. 
	Readers familiar with these topics may choose to begin at a later section and refer back to this section for notation if necessary.

	Notably, we treat \textbf{polyhedra} as closed intersection of finitely many halfspaces in $\mathbb{R}^n$, and allow unbounded polyhedra.	In our setting, a \textbf{polyhedral complex} $\mathcal{C}$ is a set of polyhedra in $\R^n$ which is closed under taking faces, such that every pair of polyhedra shares a common face (which may be the empty face).  We denote the \textbf{underlying set} of a polyhedral complex by $|\mathcal{C}|$, and take the \textbf{interior} of a polyhedron to be its interior in the relative topology induced by $\R^n$.
	A function $f:\mathcal{C}\to \mathcal{D}$ may be defined by functions with domain and codomain given by the respective underlying sets. 
	Such functions are called \textbf{piecewise linear} on  $\mathcal{C}$ if the function on $|\mathcal{C}|$ is continuous and, for each polyhedron $C \in \mathcal{C}$, $f|_C$ is affine. 
	We refer the reader to \cite{GrigsbyLindsey} and \cite{grunert}  for an additional overview of definitions in polyhedral geometry relevant to this work.

	\subsection{Some piecewise linear and polyhedral constructions}

	The terms below are defined specifically for polyhedra and polyhedral complexes. We use definitions from \cite{RS} and  \cite{Fugacci}. 
	These definitions for polyhedral complexes are motivated by similar constructions for simplicial complexes. 
	The additional generalizations to local versions of these terms are necessary for working in the polyhedral setting. We will use the local star and link of a vertex to provide combinatorial regularity which we will exploit in the algorithms we describe in \Cref{sec:relperfect} and \Cref{sec:computation}. 
	
	\begin{definition}[Cone, Cone Neighborhood, cf. \cite{RS, grunert}] 
		\label{def:cone}
		Let $p$ be a point in $\R^n$ and $A$ a set in $\mathbb{R}^n$.  We define $pA = \{t p + (1-t)a \mid t \in [0,1], a \in A  \}$ and call $pA$ a \textbf{cone} if each  point in $pA$ can be written uniquely as a linear combination of $p$ and an element of $A$.		
		A \textbf{cone neighborhood} of a point $p$ in a polyhedral complex $\mathcal{C}$ is a closed neighborhood of $p$ in $|\mathcal{C}|$ given by a cone $pA$, for a compact set $A$. 
	\end{definition}

	\begin{remark}
	 Every open neighborhood of $p$ in $|\mathcal{C}|$ contains a cone neighborhood of $p$ in $\mathcal{C}$. 
	\end{remark}

	Next, we create local versions of the following constructions which are well known in simplicial complexes. (1-2) are adapted from \cite{grunert}, (3) is adapted from \cite{Fugacci}, and (4) is, to our knowledge, new:

	\begin{definition}[Star, Local Star, Lower Star, Local Lower Star] 
		\label{def:star} \;
		
		\begin{enumerate}
			\item 	The \textbf{star}  $\st(p)$ of a point $p$ in a polyhedral complex $\mathcal{C}$ is the set of all polyhedra in $\mathcal{C}$ which contain $p$, that is, $\st(p)=\{ C \in \mathcal{C}: p \in C\}$. 	
			\item  Let $L$ be a compact set contained in $\st(p)$ such that the cone neighborhood $pL$ satisfies $pL \subset \st(p)$. The \textbf{local star of $p$ with respect to $L$},  denoted $\st_L(p)$, is given by the cells $\{C\cap pL: C \in \st(p) \}$. 
			\item If $f:\mathcal{C}\to \mathbb{R}$ is a piecewise linear function, the \textbf{lower star of p relative to $f$} is the set $\lstar(p) := \{C \in \st(p): f(x)\leq f(p) ~ \forall x \in C\}.$
			\item Finally, the \textbf{local lower star of $p$ with respect to $L$ and relative to $f$}, denoted $\lstar_L(p)$, is given by the restriction of a lower star of $p$ to the cone neighborhood $pL$: $\{C\cap pL: C \in \st(p)\}$. 
		\end{enumerate}
	
	\end{definition}

\begin{figure} \centering 
	\includegraphics[width=0.4\textwidth]{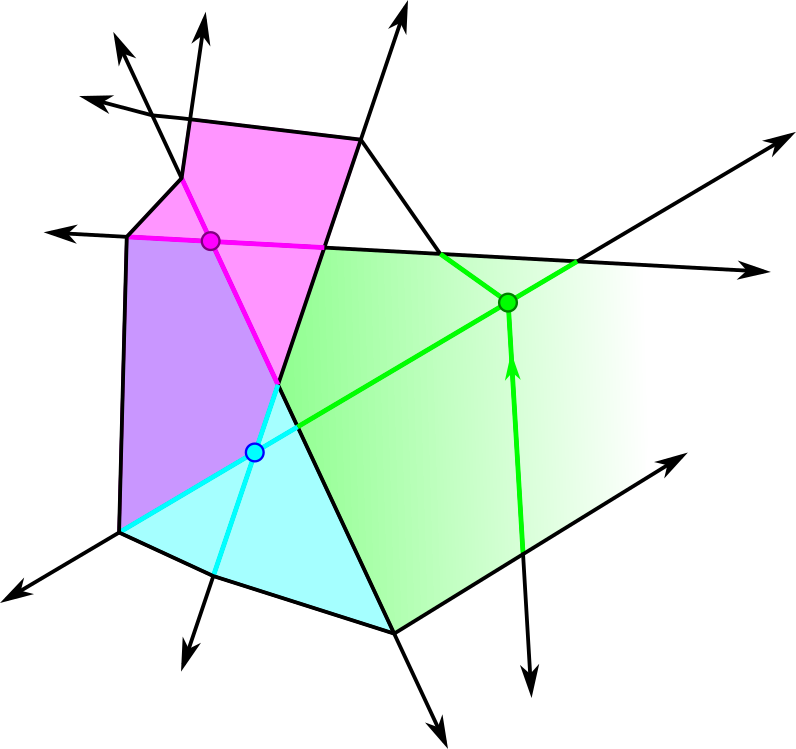} \hspace{.05\textwidth}\includegraphics[width=0.4\textwidth]{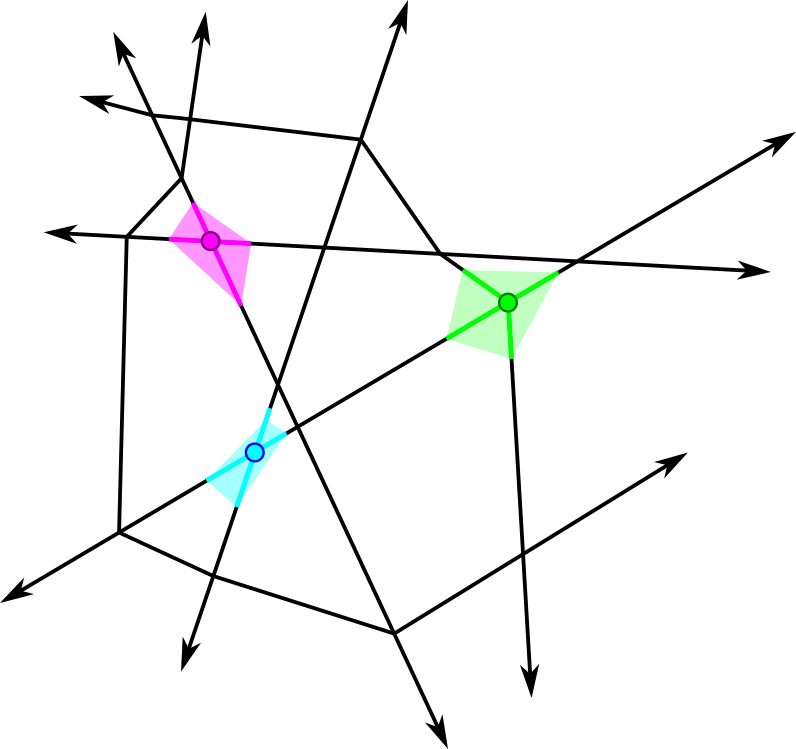} 
	
	\includegraphics[width=0.4\textwidth]{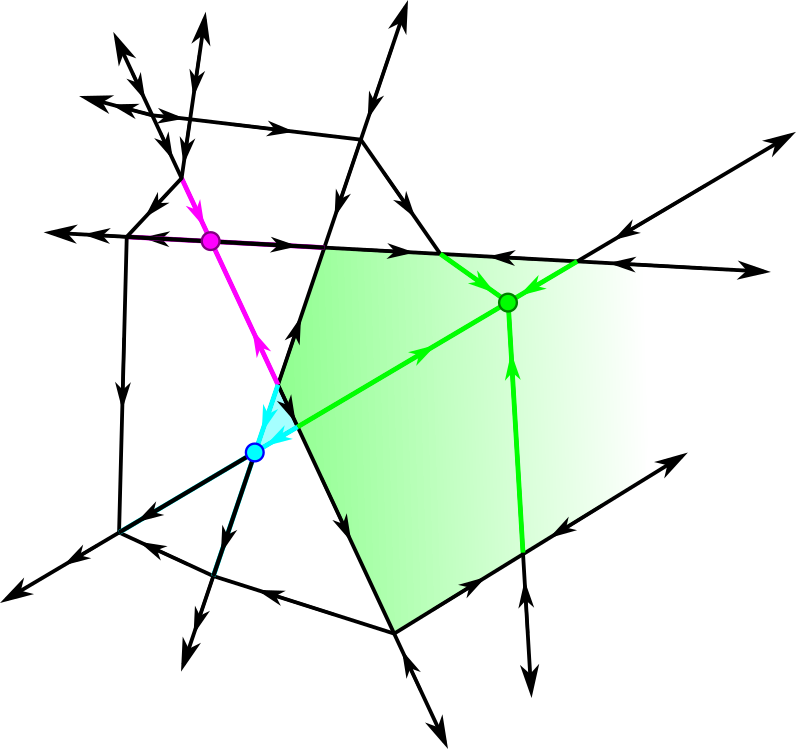} \hspace{.05\textwidth}\includegraphics[width=0.4\textwidth]{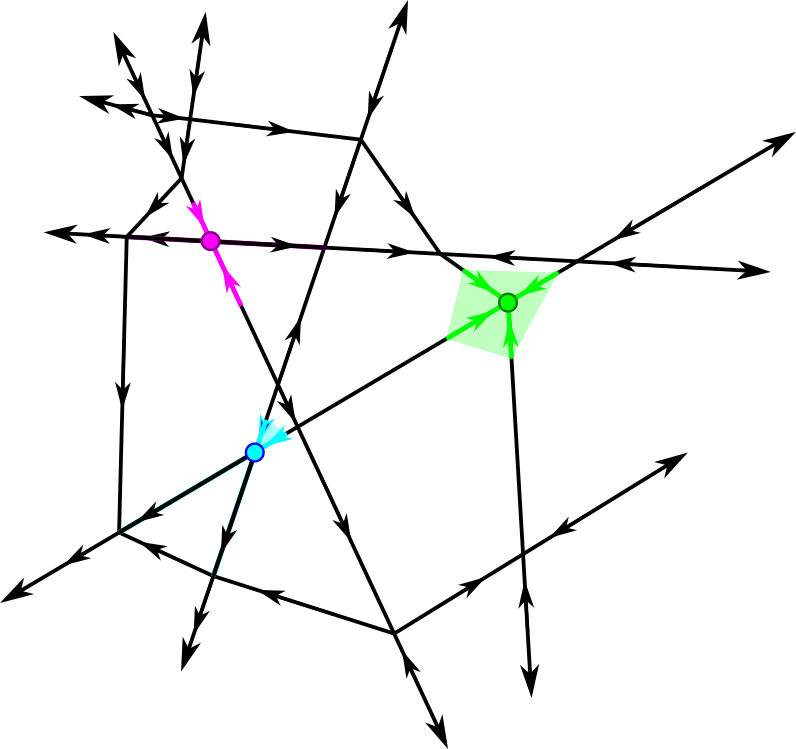} 
	
	\caption{\textbf{Upper left:} The \textit{stars} of the indicated vertices overlap.
	\textbf{Upper right:} One possible polyhedral construction of the \textit{local stars} of the indicated vertices. 
	\textbf{Bottom left:} The \textit{lower stars} of the indicated vertices, given the indicated gradient directions. Observe the lower stars are necessarily disjoint. 
	\textbf{Bottom right:} The \textit{local lower stars} of the indicated vertices have simpler combinatorial type, but the cells are in bijection with the cells in the lower star. }

	\label{fig:stardefinitions}
\end{figure}

	\begin{figure} \centering 
		\includegraphics[width=0.4\textwidth]{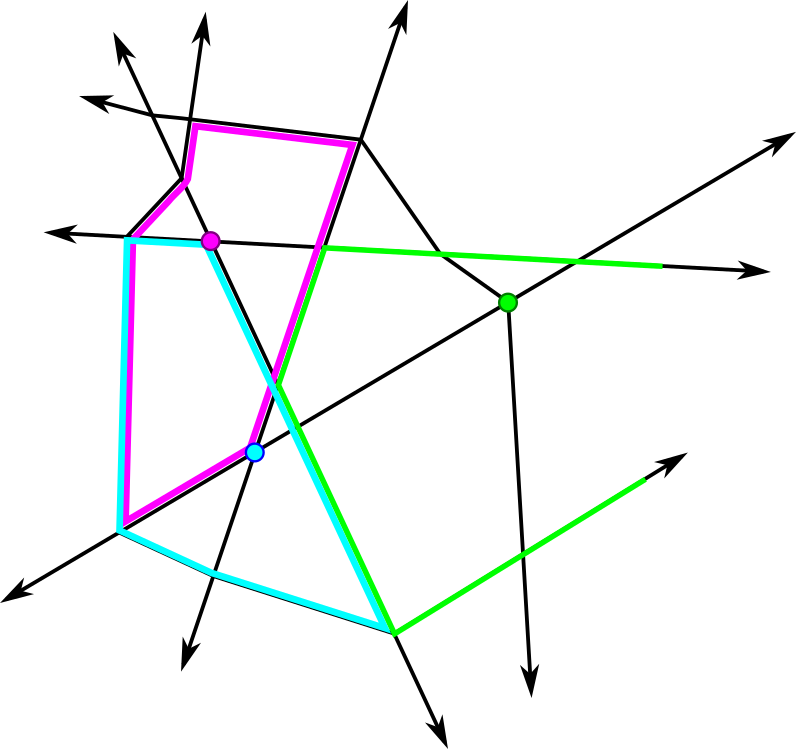} \hspace{.05\textwidth}\includegraphics[width=0.4\textwidth]{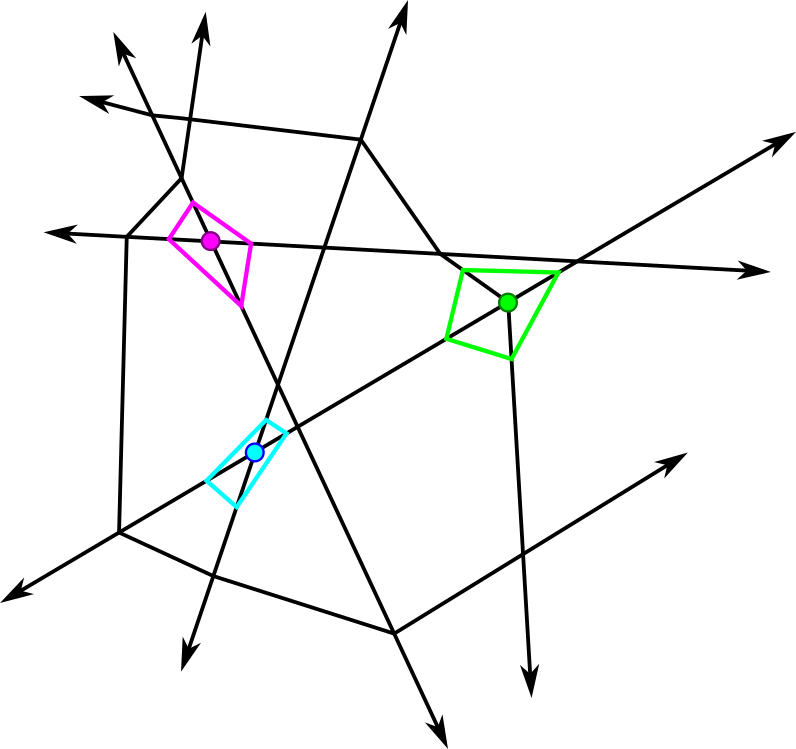} 
		
		\includegraphics[width=0.6\textwidth]{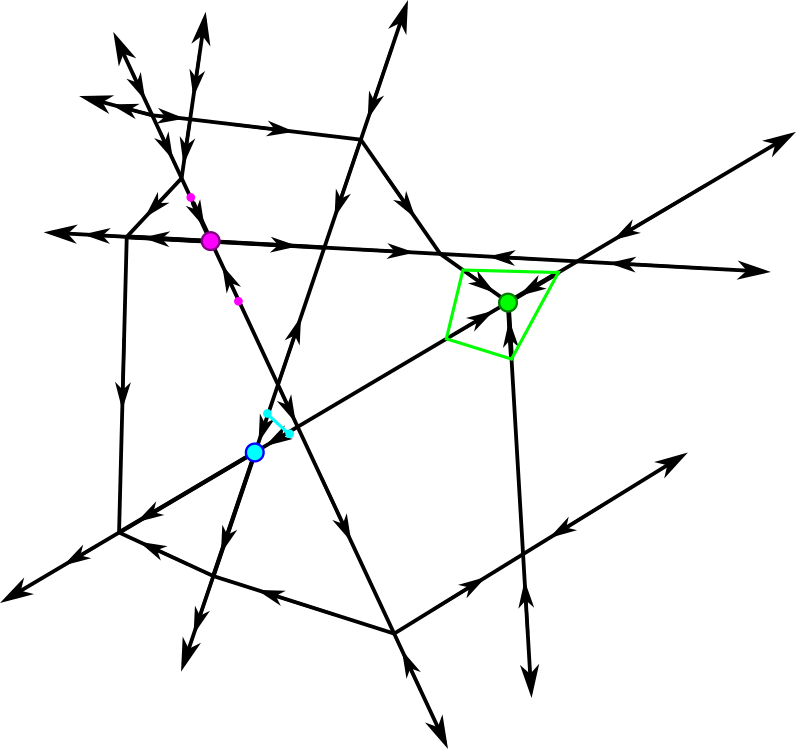} 
		
		\caption{\textbf{Upper left:} The \textit{links} of the indicated vertices overlap and have arbitrary combinatorial type. \textbf{Upper right:} The \textit{local links} of the indicated vertices. \textbf{Bottom:} The \textit{local lower links} of the indicated vertices, given the indicated $\nabla F$-orientations on edges.}
		
		\label{fig:linkdefinitions}
	\end{figure}

	Observe that for a given point $p$ and any cone neighborhood $pL$ of $p$, the local star and local lower star of $p$ have a poset structure given by containment and induced by that of the star of $p$; 
	this poset structure is independent of choice of $L$. 
	This justifies calling our construction \textit{the} local lower star of $p$ relative to $f$.

	In the simplicial context, the \textit{link} can be thought of intuitively as the boundary of the star. We introduce a local version of the link in the polyhedral setting. However, in this polyhedral setting, in contrast to the star, the link and the local link may be combinatorially distinct, even though the star and local star are not. 

	\begin{definition}[Link, Local Link, Local Lower Link] Let $p$ be a point in $|\mathcal{C}|$ for some polyhedral complex $\mathcal{C}$.
	\begin{enumerate}
		\item The \textbf{link of $p$} is the set of all faces of cells in $\st(p)$  that do not contain $p$, denoted $\lk(p)$.
		\item If $L$ is a compact set such that $pL$ is a cone neighborhood of $p$, then we call $L$ a \textbf{local link of $p$}.
		\item If $L$ is a local link of $p$ contained in $\st(p)$, then the \textbf{local lower link of $p$} is the restriction of $L$ to the lower star of $p$: $\{L \cap C: C \in \st^-(p) \}$, denoted $\llk_L(p)$.
	\end{enumerate}
	\end{definition}

	The local link and local lower link of a point $p$ in $\mathcal{C}$ have a well-defined combinatorial decomposition induced by that of the star of $p$, while the true (combinatorial) link in $\mathcal{C}$ does not, as illustrated in Figures \ref{fig:stardefinitions} and \ref{fig:linkdefinitions}.

	\subsection{ $\nabla F$-orientations of PL functions}
	
	In \Cref{sec:relperfect} we will translate from piecewise linear Morse functions on certain PL manifolds to discrete Morse gradient vector fields on a corresponding polyhedral complex. A useful intermediary, and furthermore a useful tool for visualization, is what we call the \textit{$\nabla F$-orientation} on the PL manifold's 1-skeleton.

	\begin{definition}[$\nabla F$ orientation, \cite{GrigsbyLindsey}]
		Let $\mathcal{C}$ be a polyhedral complex and $\mathcal{C}^{(1)}$ be the 1-skeleton of $\mathcal{C}$. Let $F:|\mathcal{C}|\to \mathbb{R}$ be a piecewise linear function on $\mathcal{C}$. Then the following orientation on the 1-skeleton of $\mathcal{C}$ is called the \textbf{$\nabla F$-orientation} (read ``grad-F orientation") on $\mathcal{C}^{(1)}$.
		\begin{enumerate}
			\item Orient the edges of $\mathcal{C}$ on which $F$ is nonconstant in the direction of increase of $F$.
			\item Do not assign an orientation to those edges of $\mathcal{C}$ on which $F$ is constant. 
		\end{enumerate}
	\end{definition}

\
We primarily consider the case where $F$ is only constant on vertices. For ReLU neural networks, this is not always the case, but is sufficiently common to treat as a distinguished case.

	\begin{lemma}\label{lem:gradFproperties}
		The following are properties of a $\nabla F$-orientation on a polyhedral complex $\mathcal{C}$ on which $F$ is only constant on vertices. 	
		\begin{enumerate}
			\item There are no directed cycles in the directed graph. 
			\item If $C$ is a closed, bounded polytopal cell of $\mathcal{C}$, then the $\nabla F$ orientation restricted to the boundary of the cell has a unique source and a unique sink. 
			\item If a polytope $C$ has a source $v_{min}$ and a sink $v_{max}$ induced by the $\nabla F$ orientation on its edges, then there is a hyperplane sweep in $\mathbb{R}^{n_0}$ that hits vertex $v_{min}$ first and vertex $v_{max}$ last.
		\end{enumerate}
		
	\end{lemma}
	
	\begin{proof}
		These are standard results in linear programming.  (3) is due to the fact that $F$ induces a linear projection $F\big|_C: C \to \R$ such that $F\big|_C(v_{min})$ is minimal and $F\big|_C (v_{max})$ is maximal. The preimage of each point in $\R$ is a hyperplane intersected with $C$. 
	\end{proof}
	
	\begin{remark}
		As evidenced by (3), two combinatorially equivalent polyhedral complexes might have different sets of admissible PL gradients on their edges depending on the location of their vertices. The question of how to classify all geometrically realizable polyhedral complexes and corresponding $\nabla F$-orientations by a ReLU neural network $F$ with a fixed architecture is open.
	\end{remark}

	\begin{lemma}[Cf. \cite{GLMas} Lemma 6.1]
		\label{lem:flatedge_forcing}
		If $v_1$ and $v_2$ are two vertices connected by a flat edge $e$, and $e_1$ and $e_2$ are two edges incident to $v_1$ and $v_2$ respectively which bound the same cell $C$, then the gradients on $e_1$ and $e_2$ have the same relative orientation; they must both be pointing away from $v_1$ and $v_2$ respectively, or both be pointing towards. 
	\end{lemma}
	
	The above lemma implies the following corollary immediately.
	
	\begin{corollary}
		\label{cor:codim1_forcing}
		If $C$ is a $(k-1)$ dimensional face of a $k$-dimensional cell $D$, and $F$ is a continuous affine linear function on $C$ and $D$ such that $F(C)$ is constant, then $F(D)\leq F(C)$ or $F(D)\geq F(C)$, with equality occurring only on $C$.  
	\end{corollary}
		We also find the following \textit{no-zigzags lemma} useful in understanding realizable $\nabla F$ orientations.
	
	\begin{lemma}[No-zigzags lemma] \label{lem:nozigzags} 
		Let $C$ be a $2-$cell of a polyhedral complex and $F$ an affine function on $C$. Let $e_1$ and $e_2$ be unbounded edges of $C$ with vertices $v_1$ and $v_2$. Let $e$ be an edge of $C$ connecting $v_1$ and $v_2$. If $\nabla F$ orientation on $e_1$ is towards $v_1$ and the $\nabla F$ orientation on $e_2$ is away from $v_2$ then the $\nabla F$ orientation on $e$ is from $v_1$ to $v_2$.  
	\end{lemma}
	
	\begin{proof}
		
	First, the existence of a $\nabla F$ orientation on $e_1$ and $e_2$ implies that $F$ is nonconstant on $C$. As $F$ is affine, the level sets of $F$ in $C$ are lines. 
	
	Next, $F$ is not constant on $e$, because if it were, then by \Cref{lem:flatedge_forcing} both $e_1$ and $e_2$ would be oriented in the same direction relative to $v_1$ and $v_2$. 
 	
	For the sake of contradiction, if the edge $e$ were oriented from $v_2$ to $v_1$, then consider a level set of $F$ that includes at a point $x$ on the interior of $e$. 
	The $\nabla F$-orientations on $e_1$ and $e_2$ imply that the level set containing $x$ also intersects the interiors of $e_1$ and $e_2$ in points $x_1$ and $x_2$.  The points $x_1, x_2$ and $x$, on three different faces of $C$, cannot be contained in a line because $C$ is a polyhedron, giving a contradiction.  
	\end{proof}

	\subsection{Piecewise linear Morse critical points}
	
	Smooth Morse theory describes the classical relationship between the homotopy type of sublevel sets of a smooth function $f:M\to \R$ and what is called the \textit{index} of its critical points, the number of negative eigenvalues of the Hessian \cite{milnor}. 
	In non-smooth contexts, alternative tools are needed. 	
	One such tool is piecewise linear (PL) Morse  theory. 
	No one such theory is generally accepted, but we adapt the following definition due to its 
	ease of use in this context, adapted from \cite{grunert}:

	\begin{definition}[Piecewise Linear Morse Critical Point, Regular Point, Index]\label{def:plMorse}
		Let $M$ be a combinatorial $d$-manifold and let $f: |M|\to \mathbb{R} $ be piecewise affine on cells. Let $x \in |M|$. Let $St(d)$ be the standard cross-polytope in $\mathbb{R}^d$ centered at the origin $o$ and define $f_i: St(d) \to \mathbb{R}$ by 
		
		$$f_i(x_1,...,x_d) = \sum_{k=1}^i - |x_k| + \sum_{k=i+1}^d |x_k| $$
		
		If there are combinatorially equivalent link complexes for $x$ and $o$ contained in the stars of $x$ and $o$ such that $f - f(x)$ and $f_k$ have the same signs at corresponding vertices, then $x$ is a \textbf{PL critical point of $f$ with index $i$}. 
		
		Letting $g(x_1,..,x_d)=x_1$, we call $x$ a \textbf{PL regular point of $f$} if, likewise, there is a combinatorially equivalent link complex for $x$ and $o$ contained in the stars of $x$ and $o$ such that $f-f(x)$ and $g$ have the same signs at corresponding vertices. 
		 
		If $f$ satisfies neither condition at $x$, then $x$ is called a \textbf{degenerate critical point of $f$}.
		
		Here, the \textbf{standard cross-polytope} $St(d)$ in $\mathbb{R}^d$ is the convex hull of the points $\{(0,...,\pm 1, ...,0)\}$. One natural simplicial decomposition of $St(d)$ consists of those simplices given by the convex hull of the origin, $o$, together with one vertex $v_i$ which is nonzero in the $i$th coordinate direction. This simplicial decomposition is compatible with the piecewise linear structure of $f_i$ for all $0\leq i\leq d$.	
	\end{definition}

	We say that $f$ is \textbf{PL Morse} if all vertices are regular or critical with index $i$ for some $1 \leq i \leq d$. 
	As in smooth Morse theory, the sublevel set topology of a PL function $f$ only changes at PL critical points, and if $f$ is PL Morse, the change in homotopy type at a PL critical point of index $k$ is consistent with attaching a $k$-cell. 
	Furthermore, if $f(v)=c$ and $v$ is the only PL critical point satisfying that condition, the rank of the relative homology $H_k(f_{\leq c}, f_{\leq c-\epsilon})$ is $1$ and, for $i\neq k$, $H_i(f_{\leq c}, f_{\leq c-\epsilon})$ have rank zero \cite{GLMas}. 
	
	In most cases, PL critical points only occur at vertices of the polyhedral complex, but this relies on the function $f$ being nonconstant on positive-dimensional cells. 
	In \cite{GLMas} and \cite{Masden} it is shown that all ReLU neural networks which are nonconstant on positive-dimensional cells are PL Morse.
	
	ReLU neural networks often have cells of their \textit{Canonical Polyhedral Complex} (see \Cref{def:canonicalpoly}) on which they are constant (which we call \textbf{flat cells}, in line with \cite{GLMas}), and as a result not all ReLU neural networks are PL Morse. While it is a goal of the authors to extend our construction of a discrete Morse function for networks with flat cells (see \Cref{sec:conclusion} for a brief discussion on this topic), such networks are not in the scope of the current paper.  Therefore, we restrict this paper to ReLU neural networks whose only flat cells are vertices. 
	
			
	
	
	\subsection{Discrete Morse vector fields}
	
	The strength of discrete Morse theory  is that it provides algorithmic tools for computing complete sublevel set topology \cite{Forman}. 
	 
	Such tools, to our knowledge, have not been developed in generality for PL Morse theory. 
	While the majority of discrete Morse theory has been developed in the context of simplicial complexes and CW complexes, the definitions are applicable in the context of polyhedral complexes with few changes, which we discuss at the end of this section. 
	The main difference between polyhedral complexes and cellular complexes is the presence of unbounded polyhedra, which we address in \Cref{sec:unbounded}. 
	
	We begin by reviewing standard definitions in discrete Morse theory; interested readers can find more details in \cite{Fugacci}. 

	\begin{definition}[Discrete Morse Function (\cite{Forman})]\label{def:DMF}
		Let $\mathcal{C}$ be a simplicial [cellular] complex.  A function $f:\mathcal{C}\to\R$ is a \textbf{discrete Morse function} if, for every simplex [cell] $\alpha\in \mathcal{C}$,
		
		\[
		\#\{\beta<\alpha, \textrm{dim}(\beta)=\textrm{dim}(\alpha)-1|f(\beta)\geq f(\alpha)\}\leq 1
		\]
		and
				\[
		\#\{\gamma>\alpha,\textrm{dim}(\gamma)=\textrm{dim}(\alpha)+1|f(\gamma)\leq f(\alpha)\}\leq 1,
	\]
	and at least one of the above equalities is strict.  Simplices [cells] for which both equalities are strict are called \emph{critical}.
	\end{definition}
	
	 A discrete Morse function has the property that it assigns higher values to higher dimensional simplices, except possibly with one exception locally for each simplex.  
	 For a given simplicial complex $K$, there may be many discrete Morse functions; for example, any function which assigns increasing values to simplices with increasing dimension will be discrete Morse, and all simplices will be critical.  
	 However, it is usually possible to find a more ``efficient" discrete Morse functions in the sense that complex has a smaller amount of critical simplices than overall simplices (see for example \cite{FORMAN1998,HERSH2005, Lewiner2003}).
	 
	 Sublevel sets of discrete Morse functions are subcomplexes, and provide a way of building the simplical complex in question by adding simplices in order of increasing function value.  
	 As in classical Morse theory, the homotopy type of the sublevel set can only change when a critical simplex is added, and the dimension of the critical simplex indicates that homotopy type of the sublevel differs only by the attachment of a cell of the same dimension.

	Associated to each discrete Morse function is a \emph{discrete gradient vector field}.  
	The information contained in the gradient vector field of a discrete Morse function is sufficient to determine the homotopy type of sublevel sets, and gradient vector fields have the benefit of a combinatorial formulation.  
    To define the gradient vector field of a discrete Morse function, we first introduce the definition of a general discrete vector field.
    
	\begin{definition}[Discrete Vector Field, \cite{Fugacci}]\label{def:DVF}
		Let $\mathcal{C}$ be a simplicial [cellular] complex. A \textbf{discrete vector field $V$ on $\mathcal{C}$} is a collection of pairs $\{(C,D)\}$ of simplices [cells] of $\mathcal{C}$ such that: 
		\begin{enumerate}
			\item $\dim C = \dim D- 1$ 
			\item Each simplex [cell] of $\mathcal{C}$ belongs to at most one pair in $V$.
		\end{enumerate}
	\end{definition} 
    
    The discrete gradient vector field of a discrete Morse function $f$ on $\mathcal{C}$ is the pairing that arises from $f$ as follows: 
    \begin{enumerate}
        \item If $\alpha$ is critical, then it is unpaired, 
        \item Otherwise, if there exists a face $\beta$ of $\alpha$ with $f(\beta)\geq f(\alpha)$, then $\alpha$ is paired with $\beta$, 
        \item Otherwise, there is a coface $\gamma$ of $\alpha$ with $f(\alpha)\geq f(\gamma)$, and $\alpha$ is paired with $\gamma$. 
    \end{enumerate}
    
    It is clear from  \Cref{def:DMF,def:DVF} that this collection of pairs of $\mathcal{C}$ is indeed a discrete vector field.  
    
    Often, discrete Morse functions are represented only by their gradient vector fields, as opposed to giving the function $f$ itself.  
    Therefore, it is valuable to be able to determine when a discrete vector field represents the gradient of a discrete Morse function.  
    The discrete Morse analogue to gradient flow along a discrete vector field is a $V$-path.
    
	\begin{definition}[$V$-path, \cite{Fugacci}]
		Given a discrete vector field $V$ on a  simplicial [cellular] complex $\mathcal{C}$, a \textbf{$V$-path} is a sequence of cells: 
		
		$$C_0, D_0, C_1, D_1, \ldots , C_r, D_r $$
		
		such that for each $i = 0,...,r$, $(C_i, D_i)\in V$ and $D_i>C_{i+1} \neq C_i$. 
		
	\end{definition}

	A classical way of determining whether a smooth vector field can represent a gradient is when it lacks circulation. The analogue in the discrete setting to ``lacking circulation" is ``no non-trivial closed $V$-paths".  
	The following theorem therefore indicates when a discrete vector field is the gradient of a discrete Morse function. 
	
	\begin{theorem}[\cite{Forman}, Thm. 3.5]\label{thm:forman}
		A discrete vector field $V$ on a simplicial [cellular] complex is the gradient vector field of a discrete Morse function if and only if there are no non-trivial closed $V$-paths.
	\end{theorem}

	The critical cells of a discrete gradient vector field $V$ can be thought of as tracking an analog of the topology of the polyhedral complex's sublevel sets under a particular function which would induce said discrete vector field, as seen below.

	\begin{definition}[Perfect DGVF, Relatively Perfect DGVF \cite{Fugacci}]
		\label{def:relperfect}
		
		A discrete gradient vector field $V$ on simplicial [cellular] complex $\mathcal{C}$ is called \textbf{perfect} if the number of critical cells of $V$ of dimension $k$ is equal to the rank of $H_k(|\mathcal{C}|)$ for all integers $k$. 
		
		Let $f:|\mathcal{C}|\to \mathbb{R}$ be a piecewise linear function on cells of $\mathcal{C}$. Let $\ell \in im \;f$ be restricted to the images of vertices of $\mathcal{C}$. For a given value of $\ell$, define $\ell'$ to be the greatest value of $f$ on vertices strictly less than $\ell$. 
		
		 A discrete gradient vector field $V$ on $\mathcal{C}$ is called \textbf{relatively perfect} with respect to the function $f$ if $m_i^\ell (V) = \rk H_i(|\mathcal{C}_\ell|, |\mathcal{C}_{\ell'}|) $, 
		
		$$\mathcal{C}_\ell = \{C \in \mathcal{C}\mid f_{max}(C) \leq \ell \}$$ 
		
		where $m_i^\ell(V) $ denotes the number of discrete critical $i$-simplices [cells] in $\mathcal{C}_\ell \setminus \mathcal{C}_{\ell'}$, 
	\end{definition}

	Ultimately, we hope the existence of a relatively perfect discrete gradient vector field will enable improved computation of topological features of the level and sublevel sets of a neural network $F$.
	
	\begin{figure}\centering
		\includegraphics[width=0.4\textwidth]{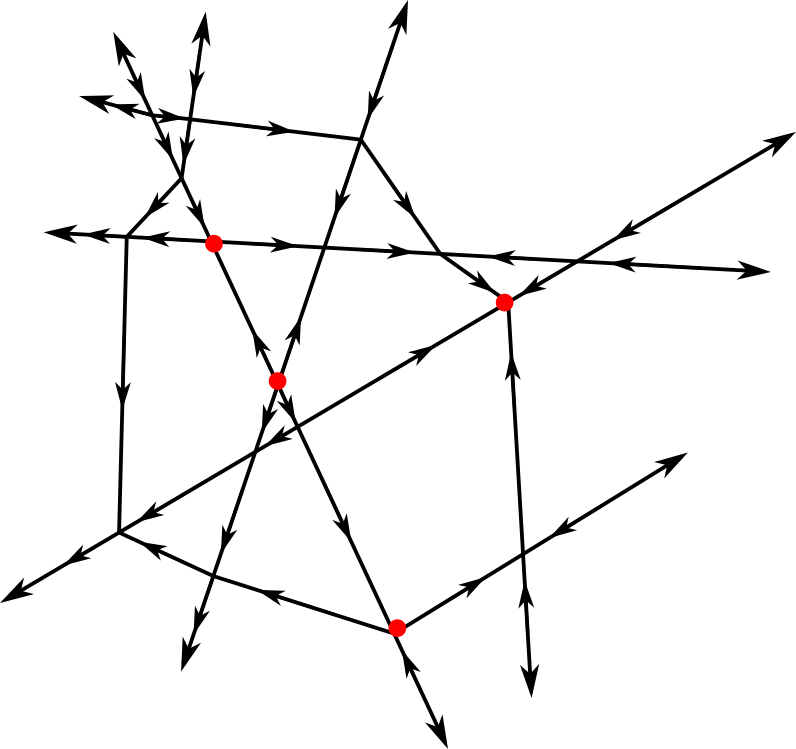} \hspace{.05\textwidth} \includegraphics[width=0.4\textwidth]{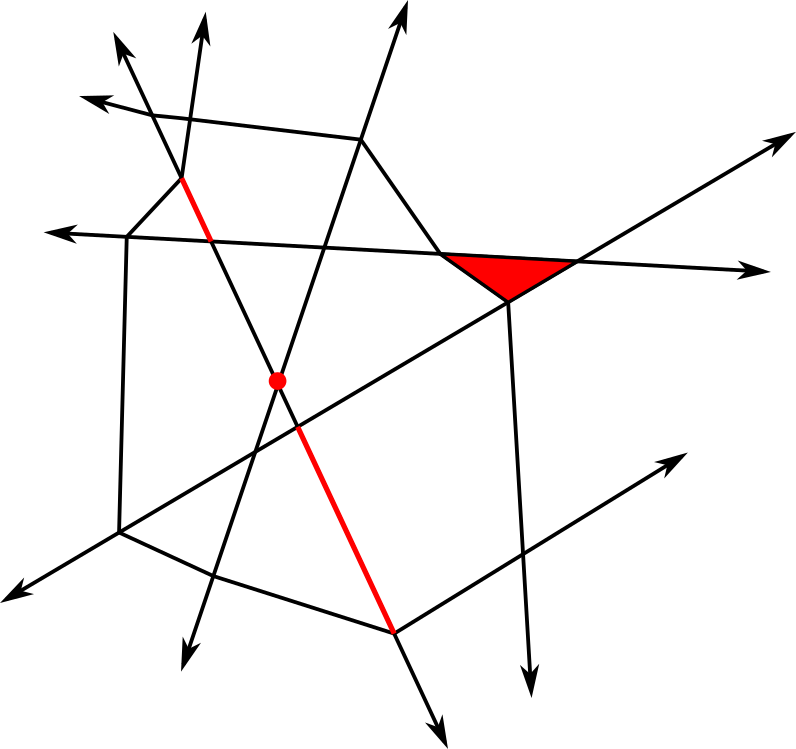}
		\caption{Left: A $\nabla F$-induced orientation on the edges of a polyhedral complex, with PL Morse critical points indicated. There are one index-zero critical point, two index-one critical points, and one index-two critical point. Right: One possible set of critical cells which would make a discrete gradient vector field on the polyhedral complex relatively perfect to $F$.}
	\end{figure}

	\subsection{Challenges in relating PL Morse and discrete Morse functions}
	
	Most known relationships beween PL Morse and discrete Morse constructions are implicit.  In \cite{Brost2013, Harker2021, grayscale}, we see that the construction of discrete gradient vector fields for cubical complexes from function data is well-studied; this is in part due to the regularity of their local combinatorics. For cubulations of compact regions of $\R^n$, for example, when simplifying persistence homology computations, there is a discrete Morse function constructed which has an implied analog of a smooth or piecewise linear function on the underlying space. 
	
	To our knowledge, there has been little work done to construct a discrete gradient vector field from a PL Morse function in a general setting. In a general setting, \cite{Lewiner} compares the PL approximations of a scalar function defined on a triangulated surface to a discrete gradient vector field built using a greedy algorithm - they find that under certain regularity conditions,  critical cells in the vector field are adjacent to critical vertices in the PL approximation. Likewise in \cite{Fugacci} an algorithm is presented for constructing a discrete gradient vector field which is relatively perfect to a PL Morse function on a simplicial complex which is a combinatorial manifold. However, due to theoretical limitations for algorithms on $n$-spheres for $n\geq 4$, the algorithm applies only to simplicial complexes of dimension $\leq 3$.  	
	
	As polyhedral complexes have fewer combinatorial restrictions on their structure than simplicial and cubical complexes, a general theory for creating discrete gradient vector fields on an arbitrary polyhedral complex from function values is likely to be similarly intractable. 
	
	Fortunately, due to the combinatorial regularity of the canonical polyhedral complex of a ReLU neural network, which we will describe in \Cref{sec:canonicalpoly}, we may follow an approach similar to that in \cite{Fugacci} to constructively obtain a relatively perfect discrete gradient vector field when the network has vertices in general position.
	
	\section{Background: The canonical polyhedral complex $\CF$}
		
	\label{sec:canonicalpoly}
	
	We now may discuss the specifics of the canonical polyhedral complex, beginning with its construction through a brief description of the combinatorial characterization which enables its topological properties to be studied.

	\subsection{Construction of $\CF$}
	
	For $n\in\mathbb{N}$, define the function $\relu:\R^n\to\R^n$ as $$\relu(x_1,\dots,x_n)=(\textrm{max}\{0,x_1\},\dots,\textrm{max}\{0,x_n\}).$$

	\begin{definition}[ReLU Neural Network; \cite{GrigsbyLindsey}, Definition 2.1,\cite{Masden}, Definition 3]\label{def:relu}
	
	Let $n_0, \dots, n_m\in \mathbb{N}$.  A \textbf{(fully-connected, feed-forward) ReLU neural network} with \textbf{architecture} $(n_0,\dots,n_m,1)$ is a collection $\mathcal{N}=\{A_i\}$ of affine maps $A_i:\R^{n_i}\to \R^{n_{i+1}}$ for $i=0,\dots,m$.  Such a collection determines a function $F_{\mathcal{N}}:\R^{n_0}\to\R$, the \textbf{associated neural network map}, given by the composite
	
	$$\R^{n_0}\xrightarrow{F_1=\relu\circ A_1}\R^{n_1}\xrightarrow{F_2=\relu\circ A_2}\dots\xrightarrow{F_m=\relu\circ A_m}\R^{n_m}\xrightarrow{G= A_{m+1}}\R.$$
	
	We say that this network has \textbf{depth} $m+1$ and \textbf{width} max$\{n_1,\dots,n_m,1\}$.  The maps $F_k$ are called the $k$th \textbf{layer maps}.  By abuse of notation, we often refer to $F_{\mathcal{N}}$ as simply $F$.    
	\end{definition}
	
The terms fully-connected and feedforward are machine learning terms which indicate that there are no restrictions on each affine function $A_i$ and that there are no expectation of identical layer maps (i.e. no recurrence), respectively. 
We will omit these terms in the rest of this paper, but keep them in the definition to disambiguate for readers with a machine learning background.

Note that $F$ is a piecewise linear function, and that $F$ defines a polyhedral decomposition of its domain, $\R^{n_0}$. 
More specifically, we may decompose $\R^{n_0}$ into a polyhedral complex by identifying the (maximal) subsets of $\R^{n_0}$ on which each $F_i$ is affine linear.   
We call this decomposition the \emph{canonical polyhedral complex of $F$}, denoted $\CF$.  
To give a precise definition of this complex, we first introduce notation concerning partial compositions of the layer maps.  
In particular, the canonical polyhedral complex can be defined using such partial composites.  
	
	\begin{definition}[\cite{Masden}, Definition 4]
	
	If $F=G\circ F_m\circ\dots\circ F_1$ is a ReLU neural network with $F:\R^{n_0}\to\R$, then we denote the composition of the first $k$ layers as $F_{(k)}$; i.e.
	$$F_{(k)}=F_k\circ\dots\circ F_1.$$
	
	We refer to $F_{(k)}$ as \textbf{ $F$ ending at the $k$th layer}.
	
	Conversely, we denote the composition of the last $m+1-k$ layers as $F^{(k)}$; i.e.
	$$F^{(k)}=G\circ F_m\circ\dots\circ F_k.$$
	
	We refer to $F^{(k)}$ as \textbf{ $F$ starting at the $k$th layer}.
	
	\end{definition}
	
Clearly, $F=F^{(k)}\circ F_{(k-1)}$.   
Furthermore, each $F_{(k)}$ has an associated natural polyhedral decomposition of $\R^{n_0}$. 
The canonical polyhedral complex will be defined as the common refinement of these iterative polyhedral decompositions. 
To formalize the polyhedral decomposition induced by $F_{(k)}$, note that each affine linear map has an associated polyhedral complex.

\begin{definition}[Notation $R^{(i)},\pi_j$;\cite{GrigsbyLindsey}, \cite{Masden}, Definition 4]
Let $A_i:\R^{n_{i-1}}\to\R^{n_i}$ be an affine function for $1\leq i\leq m$.  Denote by $R^{(i)}$ the polyhedral complex associated to the hyperplane arrangement in $\R^{n_{i-1}}$, induced by the hyperplanes given by the solution set to $H_{ij}=\{x\in\R^n:\pi_j\circ A_i(x)=0\}$, where $\pi_j$ is the projection onto the $j$th coordinate in $\R^{n_i}$.
\end{definition}

For $k>1$, $F_{(k)}$ is not affine linear, but instead is piecewise linear.  
Therefore, the solution sets associated to $F_{(k)}$ are not necessarily hyperplanes.  
However, it is still possible to use these solution sets to determine a polyhedral decomposition of the input space, for each $k$.

\begin{definition}[Node maps and bent hyperplanes, \cite{GrigsbyLindsey}, Definion 8.1, 6.1, \cite{Masden}, Definition 5, 6]

Given a ReLU neural network $F$, the \textbf{node map $F_{i,j}:\R^{n_0}\to\R$} is defined by
\[
F_{i,j}=\pi_j\circ A_i\circ F_{(i-1)}.
\]

A \textbf{bent hyperplane} of is the preimage of $0$ under a node map, that is, $F^{-1}_{i,j}(0)$ for fixed $i,j$.
\end{definition}

A bent hyperplane is generically a piecewise linear codimension $1$ submanifold of the domain (see \cite{GrigsbyLindsey} for more details). It is ``bent" in that it is a union of polyhedra, and may not be contractible or even connected. 
For each $F_{(k)}$, the associated bent hyperplanes induce a polyhdedral decompostion of $\R^{n_0}$ which we denote $\mathcal{C}(F_{(k)})$.  
The canonical polyhedral complex can then defined iteratively, by intersecting the regions of $\mathcal{C}(F_{(k)})$ with the polyhedral decomposition given at $F_{(k-1)}$.
	
	\begin{definition}[Canonical Polyhedral Complex $\CF$, \cite{Masden}, Definition 7]
	\label{def:canonicalpoly}
	Let $F:\R^{n_0}\to\R$ be a ReLU neural network with $m$ layers.  
	Define \textbf{the canonical polyhedral complex of $F$}, denoted $\CF$, as follows:
	
	\begin{enumerate}
	\item Define $\mathcal{C}(F_{(1)})$ by $R^{(1)}$.  
	\item Define $\mathcal{C}(F_{(k)})$ be defined in terms of $\mathcal{C}(F_{(k-1)})$ as the polyhedral complex consisting of the following cells:
	\[\mathcal{C}(F_{(k)})=\left\{C\cap F^{-1}_{(k-1)}(R) \ : \ C\in\mathcal{C}(F_{(k-1)}), \ R\in R^{(k)} \right\}
	\]
	
	Then $\CF$ is given by $\mathcal{C}(F_{(m)})$.
	\end{enumerate}
	
	\end{definition}	

	The above definition is the ``Forward Construction" of $\CF$ in \cite{Masden}.  Alternatively, there is an ``Backward Construction" which gives the same complex.  
	This definition originally appeared in \cite{GrigsbyLindsey}.

	For example, in special case of a neural network $F$ with architecture $(2,n,1)$, the canonical polyhedral complex $\CF$ is a decomposition of $\R^2$ by $n$ lines, which with full measure will fall in general position. 
	It is immediate that $\CF$ contains $2n$ unbounded edges, $2n$ unbounded polyhedra of dimension 2, and ${n\choose 2}$ vertices. 
	
	\begin{figure}[h]\centering
		\includegraphics[width=0.4\textwidth]{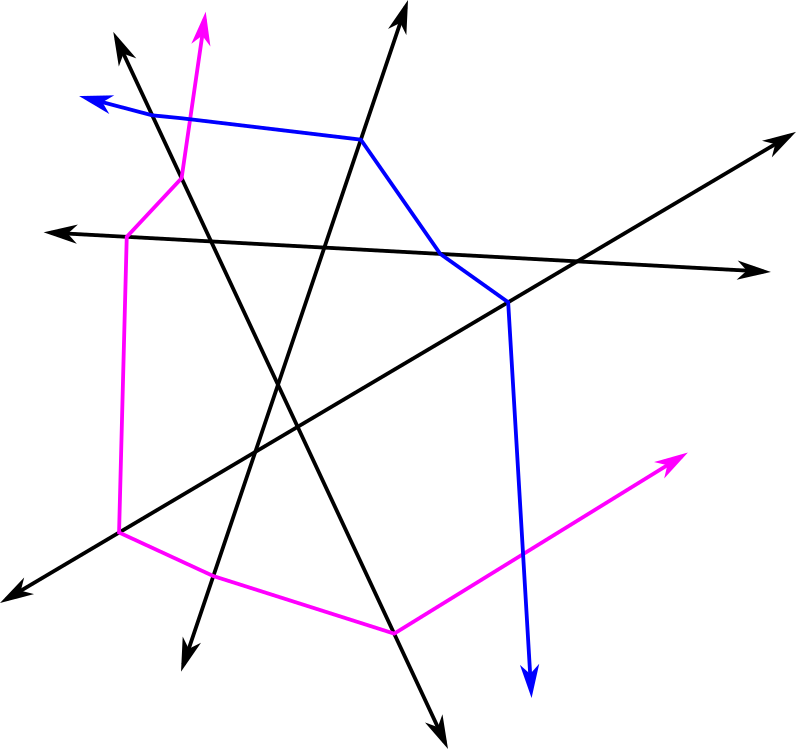}\hspace{.05\textwidth}\includegraphics[width=0.4\textwidth]{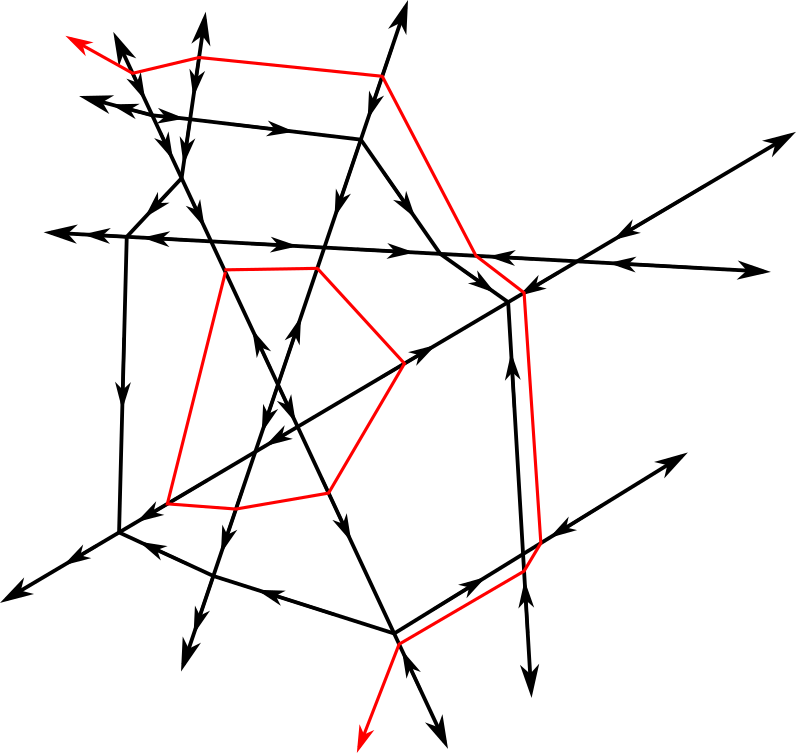}
		\caption{(Left) A portion of a canonical polyhedral complex $\mathcal{C}(F)$ given as a bent hyperplane arrangement. The two ``bent" hyperplanes which are not lines are given distinct colors. (Right) A plausible $\nabla F$ orientation on the edges of this $\mathcal{C}(F)$, and a plausible level set, marked in red.}
	\end{figure}

	\subsection{Local characterization of vertices and PL Morse critical points in $\CF$} 
	\label{sec:signsequences}
	For general piecewise linear functions on polyhedral complexes, it is not generally algorithmically decidable whether a vertex is PL critical or regular. However, the canonical polyhedral complex has combinatorial properties which make the question of PL criticality algorithmically decidable.
	
	Following \cite{Masden}, under full-measure conditions called supertransversality and genericity, each cell $C$ of $\mathcal{C}(F)$ can be labeled with a sequence in $\{-1,0,1\}^N$, where $N = \sum_{i=1}^m n_i$. This construction is not new (for example, it also appears in \cite{jordan2019provable}), but to our knowledge there is no standard reference.
	
	\begin{definition}[Sign Sequence, \cite{Masden}]
		The sign $s_{ij}(C)$ is given by the sign of $\pi_j \circ A_i\circ F_{(i-1)} (C)$, which is well-defined. The collection of all such $s_{ij}(C)$ for a specific cell $C$ is called its \textbf{sign sequence}, and is denoted $s(C)$. 
	\end{definition}

	These sign sequences encode the cellular poset of $\CF$ as follows: 
	
	\begin{lemma}[Sign sequence properties, \cite{Masden}]
		\label{lem:ssprops}
		Let $F$ be a supertransversal ReLU neural network. The following is true about any two cells $C$ and $D$ of $\CF$: 
		 		
		Define $S(C)\cdot S(D)$ by 
		
		$$(S(C)\cdot S(D))_{ij} = \begin{cases}
			S(C)_{ij} & \textrm{if $S(C)_{ij}\neq 0$} \\ 
			S(D)_{ij} & \textrm{else}
		\end{cases} $$
		
		Then $S(C)\cdot S(D)=S(E)$, where $E$ is a cell in $\CF$. Furthermore: 
		
		\begin{enumerate}
			\item $C$ is a face of $E$ (Lemma 18)
			\item $C\leq D$ if and only if $S(C)\cdot S(D)=S(D)$ (Lemma 19)
		\end{enumerate}

	\end{lemma}
	
	Finally, the cellular coboundary map in $\CF$ can be neatly described: 
	
	\begin{lemma}[Sign sequence, \cite{Masden}, Lem. 21]
		\label{lem:ss_coboundary}
	Let $F$ be a supertransversal, generic ReLU neural network. Let $C$ be a [polyhedral] cell of $\CF$. Then the cells $D$ of which $C$ is a facet are given by the set of cells with sign sequence given by $s_{ij}(D) = s_{ij}(C)$ for all $i$ and $j$ except for exactly one, a location for which $s_{ij}(C) = 0$. 
	\end{lemma}

	In other words, under supertransversality and genericity conditions each $k$-cell of $\mathcal{C}(F)$ has exactly $n_0-k$ entries of its sign sequence equal to zero, and all incident cells can be identified by replacing each zero entry with $\pm 1$. 
	Not unsurprisingly, this identifies the intersection combinatorics of cells in the local lower star of a vertex $v$ with the intersection combinatorics of the coordinate planes in $\mathbb{R}^{n_0}$. 
	
	If $v$ is a vertex of $\mathcal{C}(F)$, we can thus create a simplicial complex whose underlying set is the union of the local star and local link of $v$, and which has the combinatorics of a cross-polytope in $\mathbb{R}^n_0$. 
	The vertices of this simplicial complex are $v$ together with a point selected from each edge incident to $v$. For each $i$ from $1$ to $n_0$, replacing the $i$th zero entry with a $1$ to obtain an edge $e_i^+$, then selecting a point $v_i^+$ from that edge; likewise select $v_i^-$ from the edge $e_i^-$ obtained by replacing the $i$th zero entry with a $-1$. The $n_0$-simplices of this simplicial complex are the convex hull of the sets consisting of exactly one of $v_i^+$ and $v_i^-$ for each $i$, together with $v$. 
	
	\begin{figure}
		\centering
		\includegraphics[width=0.4\textwidth]{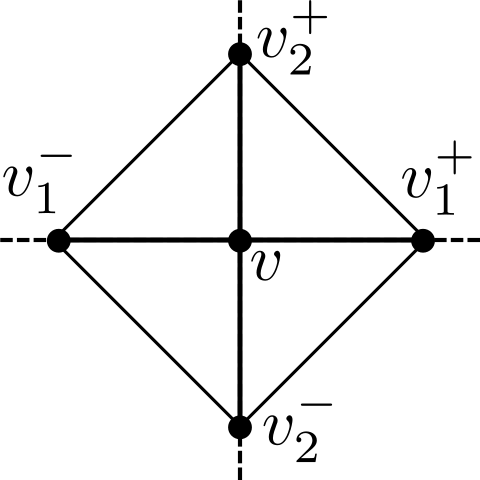}
		\caption{The union of the local star and local link of a vertex $v$ is a cross-polytope. }
		\label{fig:local_star_simplicial}
	\end{figure}

	The arguments in \cite{GLMas} and \cite{Masdendis} use this characterization and variants of \Cref{def:plMorse} to show that a vertex $v$ of $\mathcal{C}(F)$ is PL critical if and only if the edges $(v_i^-, v)$ and $(v, v_i^+)$ have opposite $\nabla F$-orientations for all $i$ , and if so, the index of the the critical point is given by the number of these pairs which are oriented towards $v$ (\cite{Masdendis}, Theorem 3.7.3). As a result, $F$ is PL Morse if and only if all edges are assigned a $\nabla F$-orientation.

	\subsection{Realizability results for PL Morse ReLU neural networks}\label{realizability}
	
	As an initial example of the usefulness of the combinatorial description given by sign sequences, we develop an exploration of some realizability results for ReLU networks by classifying all possible PL Morse ReLU neural networks on an $(n, n+1, 1)$ architecture, up to $\nabla F$-orientation (\Cref{thm:all_realizable}). These results are new in this context and our techniques are illustrative of sign sequence properties which will be used throughout the remainder of the paper. However, the main results of this paper do not rely on the results of this section.

	\begin{lemma}		\label{prop:ub_edges}
		If $F$ is a PL-Morse generic ReLU neural network with a $(n,n+1,1)$ architecture, then the $\nabla F$-orientations on unbounded edges which share the same vertex is the same, but any set of $\nabla F$-orientations on unbounded edges subject to this restriction is realizable.
	\end{lemma}

	\begin{proof}

		There is a unique generic affine hyperplane arrangment $\mathcal{A}$ with $n+1$ hyperplanes in $\R^n$, up to affine transformations. 
		It has a single bounded $n$-cell which is, in fact, an $n$-simplex, which we will call $ \Sigma$. All unbounded cells in this hyperplane arrangement share a face with this $n$-simplex. 
		
		If $F$ is PL Morse, then none of the faces of $\Sigma$ is a flat cell. This means none of the $n$-dimensional topes of $\mathcal{A}$ has the $(-1,...,-1)$ sign sequence. 
		
		There are $2^{n+1}$ possible sign sequences in $\{-1,1\}^{n+1}$, and $2^{n+1}-1$ topes in $\mathcal{A}$. 
		In particular, the $(1,1,\ldots,1)$ cell is present. Furthermore, it cannot be an unbounded cell, as the opposite unbounded cell would have the $(-1,\ldots,-1)$ sign sequence, so the $(1,\ldots,1)$ cell is $\Sigma$. Thus, the image of $F_1(\Sigma)$ is contained in the first quadrant.

		Next if $T$ is an unbounded tope of $\mathcal{A}$ with a single vertex, then it has a sign sequence with exactly one entry which is a $1$, corresponding with the axis its image is restricted to. 
		That is $F_1(T)\subset \text{span}\{e_i\}^+ $, where $i$ is the index in which the sign sequence of $T$ is positive. 
		
		All unbounded edges of $\mathcal{A}$ belong to exactly one tope $T$ of this form. An unbounded edge has exactly $n$ zero entries, and one nonzero entry, which must be $1$ otherwise the edge would be a face of the all $-1$ region. 
		The tope $T$ may be identified in sign sequence by setting all $n$ zero-entries of the edge's sign sequence to $-1$. 
		
		All edges of $T$ are unbounded and share the same vertex. A path away from the vertex of $T$ along any edge of $T$ maps to a path along $e_i$ away from the origin under $F_1$, the first layer of $F$. 
		The derivative of the restriction of $F_1$ along any edge $E$ of $T$ pointing away from the vertex of $T$ is a positive multiple of $e_i$, which we denote $D_E(F_1)=ce_i$.  
		
		The $\nabla F$-orientation on an (outward-oriented) edge $E$ of $T_i$ is given by the sign of $D_E(F_1) \cdot \vec{v} $ for $\vec{v}$ the $(n+1)$-dimensional vector giving the linear part of the affine function $G:\mathbb{R}^{n+1}\to \mathbb{R}$. 
		
		As $D_E(F_1) = c e_i$ for all edges $E$ of $T$, for a positive $c$, then the $\nabla F$-orientation on $E$ is outward if and only if the $i$th entry of $\vec{v}$ is positive. 
		
		This shows all edges $E$ of $T$ have the same $\nabla F$-orientation. 
		Since there is such a tope $T$ for each $1\leq i\leq n+1$, to induce an orientation on $T_i$ set $\vec{v}_i $ to be positive or negative, as desired. 
		This allows for all possible orientations on the unbounded edges of $\mathcal{A}$.
 	\end{proof}

  Using both \Cref{prop:ub_edges} and \Cref{lem:nozigzags}, we may further exploit properties specific to the polyhedral complex of a network $F$ with architecture $(n,n+1,1)$ to determine which vertices in $\CF$ are PL-critical, as shown in the following lemma.
  
	\begin{lemma}
		Let $F$ be an $(n,n+1,1)$ generic, PL Morse ReLU neural network. Any critical points of $F$ are index-$0$ or index-$n$, and there is at most one critical point. 
	\end{lemma}
	
	\begin{proof}
		\label{lem:max_one_crit}
		Any critical points of $F$ are vertices of $\Sigma$. Let $v$ be a vertex of $\Sigma$ and suppose it is a critical point. 
		Let $T$ be the unique unbounded tope of $\mathcal A$ whose only vertex is $v$. Note that $T$ has $n$ (unbounded) edges, and the $\nabla$-F orientation on these edges of $T$ either are all towards $v$ or are all away from $v$, as seen in  \Cref{prop:ub_edges}. 
		As each edge of $T$ is opposite $v$ from an edge of $\Sigma$, in order for $v$ to be critical all edges on $\Sigma$ containing $v$ must be oriented in the opposite direction from their paired edges on $T$ \cite{GLMas}. 
		As a result, the edges in $\mathcal{A}$ incident to $v$ are either all oriented towards $v$ or all oriented away from $v$; that is, $v$ is either index-$0$ or index-$n$.
		
		To see that there is at most one critical point, without loss of generality, assume $v$ is critical of index $0$. 
		Then all edges incident to $v$ are oriented away from $v$. If $w$ is another vertex of $\Sigma$, then there is an edge $e$ connecting $v$ to $w$, and of course $e$ must be oriented away from $v$. 
		There is a unique unbounded cell $U$ in $\mathcal{A}$ which contains $v$ and $w$ and no other vertices; it also contains $e$. The unbounded edges of $U$ which contain $v$ must be oriented away from $v$.  
		By the no-zigzags lemma \ref{lem:nozigzags} the unbounded edges of $U$ which contain $w$ must be oriented away from $w$ (as each unbounded edge of $U$ containing $w$ shares a $2$-cell with an unbounded edge of $U$ containing $v$). By  \Cref{prop:ub_edges} all unbounded edges pointing away from $w$ are oriented away from $w$. In particular, the edge opposite $e$ is oriented away from $w$. 
		Because $e$ is oriented towards $w$ and its opposite edge is oriented away from $w$, we conclude $w$ is a PL-regular point. 
	\end{proof}

	\begin{theorem}\label{thm:all_realizable}
		Let $F$ be an $(n,n+1,1)$ generic, PL Morse ReLU neural network. Then the decision boundary of $F$ is empty, has the homotopy type of a point, or has the homotopy type of an $(n-1)$-sphere. 
	\end{theorem}
	
	\begin{proof}
		Following the previous lemma \ref{lem:max_one_crit}, we will follow \cite{GLMas}, and conclude the topology of the decision boundary must be one of these three options. 
	\end{proof}

	We now classify all possible PL Morse $(2,3,1)$ neural networks.  
	
\begin{figure}[h] \centering
		\begin{subfigure}{0.45\textwidth}
			\includegraphics[width=0.7\textwidth]{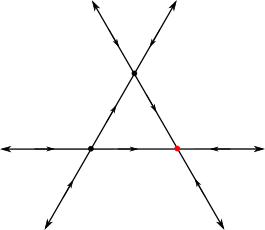} 
			\caption{$(2,3,1)$ neural network with PL-critical vertex of index $2$ marked in red.}
			\label{fig:gradF1}
		\end{subfigure}
		\hfill
		\begin{subfigure}{0.45\textwidth}
			\includegraphics[width=0.7\textwidth]{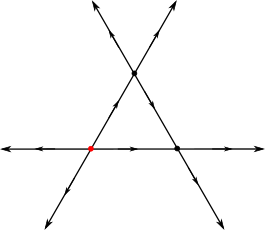} 
			\caption{$(2,3,1)$ neural network with PL-critical vertex of index $0$ marked in red.}
			\label{fig:gradF4}
		\end{subfigure}
		\begin{subfigure}{0.45\textwidth}
			\includegraphics[width=0.7\textwidth]{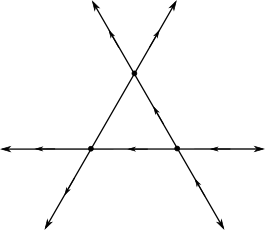} 
			\caption{$(2,3,1)$ neural network with PL-regular vertices; the unmarked edge can have either possible orientation.}
			\label{fig:gradF2}
		\end{subfigure}
		\hfill
		\begin{subfigure}{0.45\textwidth}
			\includegraphics[width=0.7\textwidth]{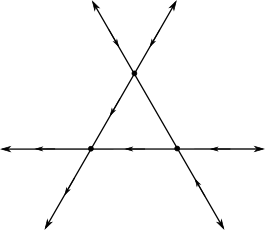} 
			\caption{$(2,3,1)$ neural network with PL-regular vertices; the unmarked edge can have either possible orientation.}
			\label{fig:gradF3}
		\end{subfigure}
		\caption{All possible $\nabla F$-orientations for a generic, supertransversal, PL Morse $(2,3,1)$ neural network.}
		\label{fig:gradF}
	\end{figure}

	\begin{corollary} \label{cor:all_231}

		The $\nabla F$-orientations depicted in \Cref{fig:gradF} are the only possible $\nabla F$-orientations on a generic, supertransversal, PL Morse $(2,3,1)$ ReLU neural network, up to combinatorial equivalence. 

	\end{corollary}

	\begin{proof}
		
		Let $F$ be such a network, and denote as $\Sigma$ the unique bounded $2$-cell of $\CF$. Following \Cref{prop:ub_edges}, there are 4 possible scenarios for the $\nabla F$-orientations on the unbounded edges of $\CF$:  
		\begin{enumerate}
			\item all unbounded edges are oriented towards $\Sigma$ (\Cref{fig:gradF1}),
			\item all unbounded edges are oriented away from $\Sigma$ (\Cref{fig:gradF4}),
			\item the unbounded edges of exactly one vertex of $\Sigma$ are oriented towards $\Sigma$ and all other unbounded edges are oriented away from $\sigma$ (\Cref{fig:gradF2}), or
			\item the unbounded edges of exactly one vertex of $\Sigma$ are oriented away from $\Sigma$ and all other unbounded edges are oriented towards $\sigma$(\Cref{fig:gradF3}).
		\end{enumerate}
		
		If not all unbounded edges have the same orientation with respect to $\Sigma$, then \Cref{lem:nozigzags} determines the orientation of two of the three edges of $\Sigma$.  Moreover, the orientations of these bounded edges ensure that none of the vertices in $\sigma$ can be PL-critical.

		If all unbounded edges have the same orientation with respect to $\Sigma$, then \Cref{lem:gradFproperties} ensures that $\nabla F$ orientation on the edges of $\Sigma$ do not generate a cycle.  Therefore, there must be a vertex $v$ of $\Sigma$ for which the $\nabla F$ orientation of each bounded edge adjacent to $v$ is towards $v$, and a vertex $w$ of $\Sigma$ for which the $\nabla F$ orientation of each bounded edge adjacent to $w$ is away from $w$. 
		
		If it is the case that all unbounded edges are oriented towards $\Sigma$, then $v$ is a PL-critical vertex of index $2$, and all other vertices are PL-regular.  
		
		If it is the case that all unbounded edges are oriented away $\Sigma$, then $w$ is a PL-critical vertex of index $2$, and all other vertices are PL-regular.  
 		\end{proof}

	Given a  canonical polyhedral complex (arising from a ReLU neural network), not all functions on that canonical polyhedral complex are realizable as ReLU neural network functions. From \cite{GLMas} we see that that a ReLU neural network of the form $(n, m, 1)$ generally has at most one $n$-cell on which it is constant, for example, but other limitations exist as well.

	\begin{figure}[h]\centering 
		\includegraphics[width=.55\textwidth]{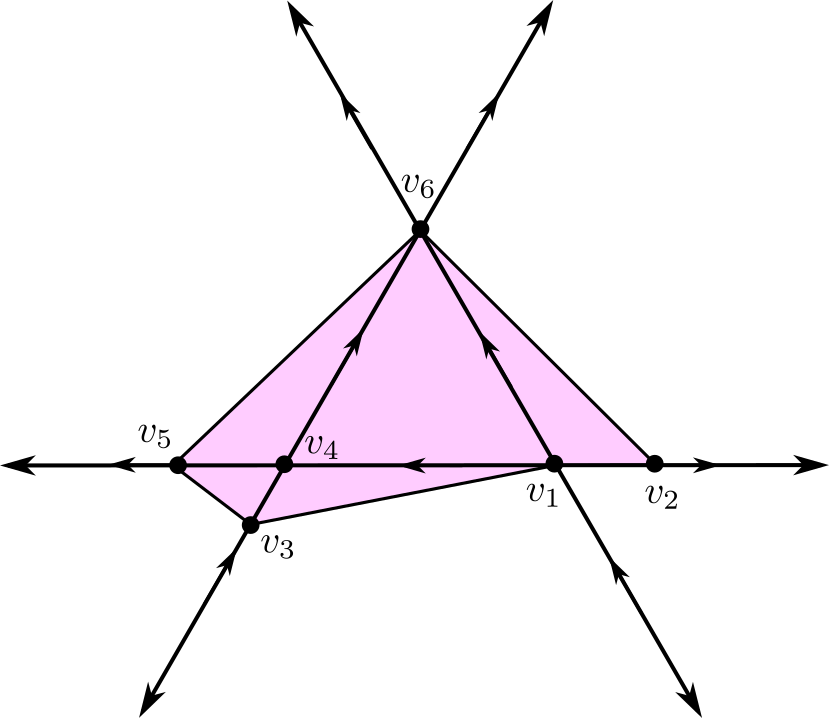}
		\caption{A $\nabla F$ orientation on the canonical polyhedral complex consisting of a generic hyperplane arrangement of 3 hyperplanes in $\R^2$ which is realizable as a PL function, but not as a ReLU neural network of with architecture $(2,3,1)$. Selecting distinct values on $v_1,...v_6$ determines a $\nabla F$ orientation on the pink polyhedral complex.}
		\label{fig:PL_not_ReLU}
	\end{figure}
	
	In fact, \Cref{fig:PL_not_ReLU} depicts a $\nabla F$ orientation on the same 3 hyperplanes which is realizable if $F$ is a PL function, but not if $F$ is a ReLU neural network of architecture $(2,3,1)$. That it cannot be a ReLU network follows immediately from \Cref{cor:all_231}. That it is realizable can be seen by assigning $F(v_i)=i$ and observing that $F$ on each simplex in the highlighted simplicial complex determines the $\nabla F$ orientation on the edges of the polyhedral cell which contains it. This results in the $\nabla F$ orientations pictured.

	\section{Relatively Perfect Discrete Gradient Vector Fields for ReLU Networks} 
	
	\label{sec:relperfect}
	
	Because our setup allows us to identify PL critical vertices, we further leverage this information to constructively establish the existence of a discrete gradient vector field on the cells of $\mathcal{C}(F)$ which are bounded above in $F$. Moreover, this discrete gradient vector field has the property that critical cells are in bijection with PL critical vertices in a way which respects function values; this is a technical property introduced in \cite{Fugacci} called \textit{relative perfectness} (\Cref{def:relperfect}). Ideally, the existence and algorithmic constructability of this discrete gradient vector field will enable faster computational measurements of the topology of ReLU neural networks' decision regions. 
	
	In this section, we follow a similar construction to that in \cite{Fugacci}, where they establish an algorithm for finding a discrete gradient vector field which is relatively perfect to any given PL Morse function on a simplicial combinatorial manifold with dimension $d\leq 3$. 
	In a similar vein, we construct a discrete gradient vector field by considering the lower stars of individual vertices of $\mathcal{C}(F)$. 
	
	Some of the key differences in our algorithm are that (a) we are not dimensionally restricted, (b) $\mathcal{C}(F)$ is not generally a simplicial complex, and (c) we do not rely on nonconstructive existence theorems to assign local gradient vector fields. 
	Instead, we exploit specific combinatorics of $\mathcal{C}(F)$ (as given in \Cref{sec:signsequences}) to constructively produce the desired local pairings. 
	To our knowledge, constructions establishing discrete gradient vector fields on polyhedral complexes with associated PL Morse functions are relatively unexplored.

        \subsection{Discrete Morse theory on unbounded polyhedral complexes}
        \label{sec:unbounded}
        
       Before we introduce our construction, we must justify why it is reasonable to use the tools of discrete Morse theory on the polyhedral complex $\CF$, which is not formally a cellular complex due to the presence of unbound cells.  In fact, we will construct a discrete gradient vector field on an associated $CW$-complex $\CF^-_*$ with \textit{no} unbounded cells.

        \begin{definition}[$\CF^-, \CF^-_*$]
        	The \textbf{Complete Lower Star Complex relative to $F$} is the subcomplex of $\CF$ containing all cells which are bounded above, i.e.
        	
        	 \[\CF^-:=\{C\in\CF  |\  \exists\  r\in\R \textrm{ such that } C\subset F^{-1}(-\infty,r])\}.
        	\]
        	
        	The one-point compactification of $\CF^-$, which we call the \emph{ one-point compactified complete lower star complex} relative to $F$ is denoted $\CF^-_*$. 
        	The distinguished point $\{*\}$ is formally assigned a function value $F(*) = -\infty$.
        \end{definition}
        
        \begin{remark}
        	$\CF^-$ is indeed a subcomplex of $\CF$, as if $C$ is a polyhedron satisfying $F(C)\leq r$, then this is true of all of its faces as well. Furthermore, even though $\CF^-$ contains unbounded cells, $\CF^-_*$ is a regular $CW$-complex. 
        \end{remark}

        In \Cref{thm:relativelyperfectDGVF} we will identify a a discrete Morse function on $\CF^-_*$ with a single connected component appearing at $-\infty$. 
        This new discrete Morse function will be relatively perfect to the PL function on the subset of $S^{n_0}$ given by $\CF^-_*$.   
	That this construction does not capture gradient information on cells whose image in $F$ is not bounded above is not a large problem: the homotopy type of the sublevel sets of $\mathcal{C}(F)$ only changes at $F(v)$, for $v$ a vertex of $\mathcal{C}(F)$. 
    
		\subsection{Networks which are injective on vertices}
				  
		We start our construction locally, by showing that, for a given vertex in $\CF$, it is always possible to  generate pairings in the local lower star of $v$ which reflect the PL-criticality of $v$.
		
		\begin{lemma} \label{lem:local-acyclic-pairing}
			Let $F$ be a fully-connected, feedforward ReLU neural network which is injective on vertices of $\mathcal{C}(F)$. 
			Then for each vertex $v$ in $\CF$, there is a pairing in the local lower star of $v$ relative to $F$ satisfying exactly one of the two following conditions:
			
			\begin{enumerate}
				\item 			If $v$ is a vertex of $\mathcal{C}(F)$ that is PL regular, then there exists a choice $V$ of complete acyclic pairing of cells in the local lower star of $v$ relative to $F$.
				\item 			If $v$ is a vertex of $\mathcal{C}(F)$ that is PL critical of index $k$, then there exists a choice $V$ of acyclic pairings of cells in the local lower star of $v$ which leaves exactly one $k$-cell unpaired. 
			\end{enumerate}
			
			Furthermore, these pairings can be constructed algorithmically. 
		\end{lemma}
		
	\begin{proof}
		As $F$ is injective on vertices, all edges have a $\nabla F$-orientation and $F$ is PL Morse.	For any vertex $v\in \mathcal{C}(F)$, following \Cref{sec:signsequences} there exists a compact set $L\subset \st(v)$ for which $\st_L(v)$ is a $n_0$-dimensional cross-polytope, and  the cells of this cross polytope are in one to one correspondence with the cells in $\st(v)$.  Through this correspondence, a pairing on the cells of $\st_L(v)$ may be used to induce a pairing on the cells of $\st(v)$.  This correspondence restricts to a correspondence between 
        $\st_L^-(v)$ and $\st^-(v)$.
					
	\textbf{If $v$ is regular (\Cref{fig:lowerstarpairing}).}   Let $v$ be a PL regular point. Then  there is at least on pair of edges $(e_i^+, e_i^-)$ for which $e_i^+\notin\st_L^-(v)$ and $e_i^-\in\st_L^-(v)$ (or vice versa) (due to \cite{Masdendis}, Theorem 3.7.3). Assume without loss of generality it is $e_i^-$ (as local relabeling does not change the combinatorics).

	\begin{figure}[h]\centering 
		\includegraphics[height=0.4\textwidth]{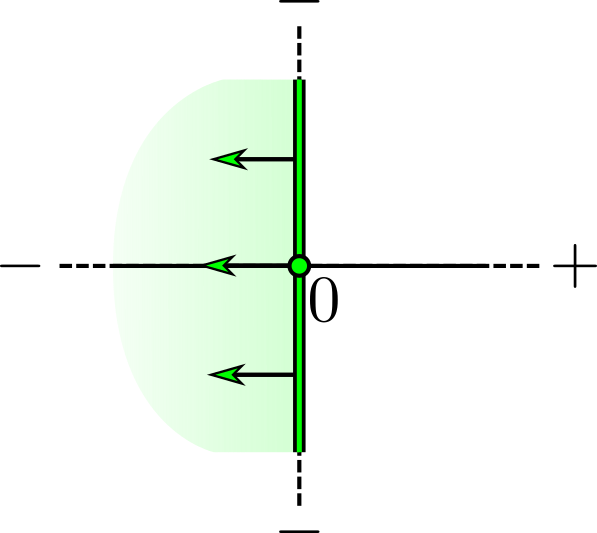}\hspace{0.1\textwidth} \includegraphics[height=0.4\textwidth]{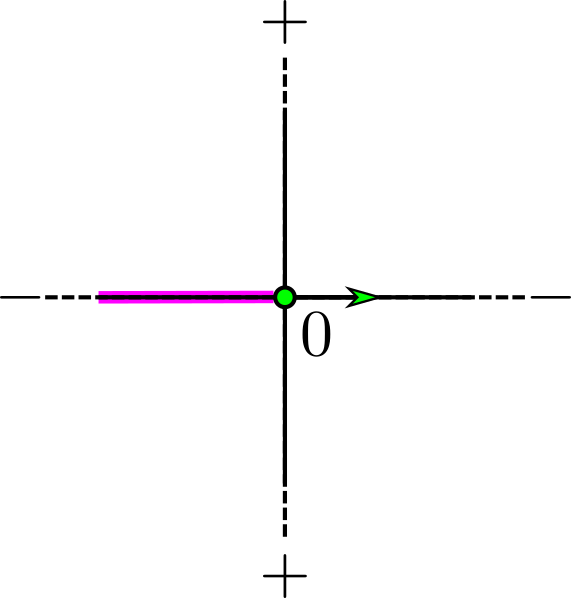}
		\caption{Left. A pairing on the lower star of a a PL regular point. Right. A pairing on the lower star of an index 1 PL critical point.  }
		\label{fig:lowerstarpairing}
	\end{figure}

 Let $v^-_i$ denote the bounding vertex of $e^-_i$ in $\st_L(v)$ which is not $v$.  We may view $\st_L^-(v)$ as the cone of $v^-_i$ over $T(v)\cap\st_L^-(v)$, where $T(v)$ is the $(n_0-1)$-dimensional cross-polytope obtained when the antipodal edges $e^+_{i^*}$ and $e^-_{i^*}$ are removed from $\st_L(v)$:  i.e.,
\[
	T(v):=\{C\in \st_L(v) | e^+_{i},e^-_{i}\notin C\}.
\] 

 We may obtain a discrete vector field $V$ on $\st_L^-(v)$ by pairing a $d$-dimensional cell $\sigma\in T(v)\cap \st_L^-(v)$ with the corresponding $d+1$ dimensional cell $v^-_{i^*}\sigma$.  This gives a complete pairing of all cells in $\st_L^-(v)$.   
 
 To see that $V$ is acylic note that, by construction, any path within $V$ contains exactly one pair. For any choice of $(C_0,D_0)$ as an initial pair in a path, we have $D_0=v^-_{i} C_0$.  Any choice of $C_1<D_0$ and $C_1\neq C_0$ has the property $C_1$ contains $v^-_{i^*}$, so that $C_1$ paired with a codimension one face in $V$.  This ensures that any path terminates after one pair.

		\begin{figure}[h] \centering
			 \label{fig:pairing_order}
			\includegraphics[width=0.4\textwidth]{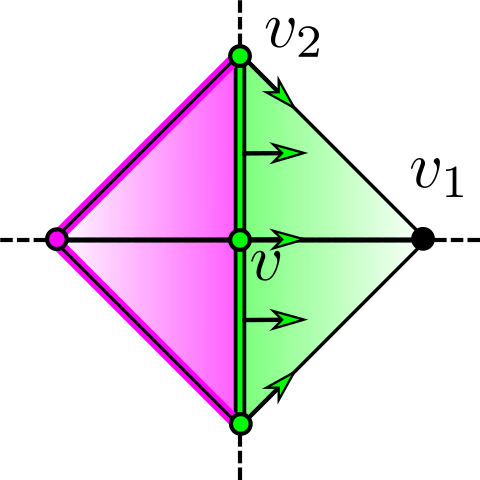} \hspace{.15\textwidth} 
			\includegraphics[width=0.4\textwidth]{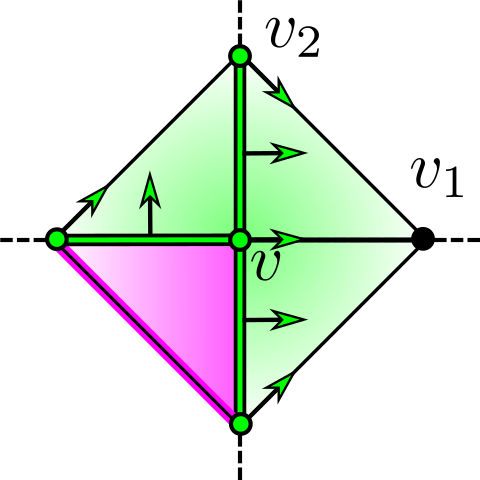}
			\caption{The algorithm constructing the pairing on $\lstar_L(v) \cup L$ for an index-$2$ critical vertex $v$. This pairing restricts to cells in $\st^-(v)$, leaving only one two-cell unpaired. }
		\end{figure}

		\textbf{If $v$ is critical} (\Cref{fig:pairing_order}).
		Now suppose $v$ is critical of index $k$. 
		Let $L$ be a local lower link of $v$. 
		Then $\lstar_L(v)$ has the combinatorics of the coordinate axes of $\mathbb{R}^k$, and is the cone of $v$ with $L$.
		 Together,  $\lstar_L(v) \cup L$ form a $k$-dimensional simplicial complex, call it $\mathcal{N}$, whose underlying set is a $k$-dimensional cross-polytope. 
		
		 \textit{Selection step.} 
		 For each coordinate direction $1 \leq i \leq k$, there is a pair of edges $e_i^+$ and $e_i^-$ in $\st^-(v)$ whose sign sequences differ by opposite signs in a single entry. 
		 The choice of labeling of $e_i^+$ determines $e_i^-$. 
		 These edges intersect $\lk^-(v)$ on opposite sides of $v$ in vertices $v_i^+, v_i^-$.
		 For each $1\leq i\leq k$, let $v_i$ be a vertex selected from $\{v_i^+, v_i^-\}$ (distinguishing a single quadrant of $\st^-(v)$).
		  Observe that each $v_i$ has an opposite vertex given by the other vertex in this set, which we will call $v_i^{op}$. 
		  Once the choice of ordering coordinate directions and choosing a positive direction in each coordinate direction is complete, no other choices need to be made.

		\textit{Pairing construction algorithm}.
		 We can then recursively assign the following pairings on $\mathcal{N}$. 
		 Recursively for $1\leq i < k$, for each $C$ in $\lk(v_i)$ (relative to $\mathcal{N}$), add the pairing $(C, Cv_i)$ if $C$ has not yet been paired. These pairings restrict to a pairing on $\st^-(v)$.  
		
		\textit{This operation constructs a discrete vector field. } 
		Namely, it does not try to pair any cell twice. 
		Observe that at step $i$ if $C$ in $\st^-(v)$ is in the $\lk_{\mathcal{N}}(v_i)$ and $C$ has not yet been paired, then we claim $Cv_i$ has also not yet been paired.
		This is because  each step adds cells which comprise the union of the star and the link of a single vertex in $\mathcal{N}$, a simplicial complex; this union is always a simplicial complex. 
		If a cell $C$ has not been added to this union, then no cell it is the boundary of could have been added to the union either. 
		Thus each cell $C$ is included in at most one pair of cells in the pairing, as is needed.

		\textit{All cells in $\lstar_L(v)$ get paired except one, of dimension $k$. }
		By construction, the union of the cells of the resulting pairing is the union of the stars and links of the vertices $v_1,...,v_n$ in $\mathcal{N}$. 
		Since $\st_{\mathcal{N}}(v_i^+)\cup\lk_{\mathcal{N}}(v_i^+)$ is precisely all cells in $\mathcal{N}$ which do not contain $v_i^-$ (and vice versa), we see that the only cell in $\lstar_L(v)$ which is not paired is the unique cell which contains none of the $v_i$ and is not in the stars of any of the $v_i$. This is the interior of the simplex given by $\{v\}\cup \{v_i^{opp}\}_{1\leq i\leq k}$, which is a $k$-simplex, as desired.
		
		\textit{The resulting construction is acyclic.}  
		Let $(C_1, D_1)$ be a pair in $V$ arising from this algorithm. Then by construction $D_1 = C_1 v_i$ for some $i$ satisfying $1 \leq i \leq k$. 
		
		Now we consider the possibilities for the next pair in the $V$-path, $(C_2, D_2)$. We know $C_2$ must be a codimension 1 element of the boundary of $D_1$ that is not $C_1$. It also must be paired with a higher-dimensional coface. 
		
		Let $C$ be a candidate for $C_2$, that is, an arbitrary codimension 1 element of the boundary of $D_1$ that is not $C_1$. We observe that, as a boundary element of $D$ which is not $C_1$, 
		
		$$C = B v_i$$ 
		
		where $B$ is a codimension-one boundary element of $C_1$. 
		
		Consider the step at which $C$ was added to the pairing. 
		
		If $C$ had not yet been added to the pairing by step $i$, then $C$ contains $v_j^{op}$ for all $j<i$ (as otherwise $C$ would be in the star or link of $v_j$).  Thus $B$ also contains $v_j^{op}$ for all $j<i$ as well, and also must not have been paired by step $i$. If $B$ was not a member of a pair by step $i$, then $(B, Bv_i = C)$ is a $V$-pair, and $C$ is added to the pairing at step $i$. Thus, $C$ is not the first element of a pairing and cannot be $C_2$ in a $V$-path.

		We conclude that if $(C_1, D_1), (C_2, D_2)$ is a $V$-path within the local lower star of $v$, then $C_2$ was paired at a strictly earlier step than $C_1$. As a result, $V$ contains no cycles.

		\end{proof}
		
		\begin{remark}
			Observe that the only choice made was in the \textit{Selection Step.} 
			To make the \textit{Selection Step} deterministic and dependent on the values of $F$, we observe that $v_i$ can ``morally" be selected to be on those edges whose opposing vertex has the lowest value in $\mathcal{C}(F)$, or if there is no opposing vertex, whose unbounded edge has the steepest directional derivative. 
			However, the same result follows regardless of whether we made the ``moral" choice. In fact, as we discuss in \Cref{sec:computation}, it is potentially more computationally convenient to be amoral (in this sense). 
		\end{remark}

		We now are able to ``stitch together" the local pairings on the lower stars of the vertices for a global discrete vector field which satisfies the desired properties, following a similar approach as in \cite{Fugacci}.

		\begin{theorem}\label{thm:relativelyperfectDGVF}
			Let $v$ be a vertex in $\CF$ and let $V(v)$ be an acyclic discrete vector field on $\st^-(v)$ obtained by the construction in \Cref{lem:local-acyclic-pairing}. Then $V = \bigcup_{v\in\CF} V(v)$ is a relatively perfect discrete gradient vector field to $F$ on $\mathcal{C}(F)^-_*$, with $\{*\}$ a critical $0$-cell.
		\end{theorem}

		\begin{proof} 
			
			The lower stars of each of the vertices of $\mathcal{C}(F)$ are disjoint, and the union of all the lower stars of all the vertices of $\mathcal{C}(F)$ is $\mathcal{C}(F)^-$.  
			
			\textit{$V$ is a valid discrete gradient vector field. }
			Suppose we have a $V$-path $(C_1, D_1),(C_2, D_2),..., (C_n, D_n)$. 
			If this path consists entirely of cells in the lower star of some fixed vertex $v$ then by \Cref{lem:local-acyclic-pairing} it is acyclic. 
			
			Otherwise, some pair $(C_i, D_i), (C_{i+1},D_{i+1})$ satisfies the condition that $C_i$ is in the lower star of a different vertex than $C_{i+1}$; call these $v_i$ and $v_{i+1}$ respectively. 
			(Observe we can make this statement because each $V$-pair is contained in the lower star of the same vertex, by construction). 
			As $D_i$ is also in the lower star of $v_i$, $F(D_i)$ is bounded above by $F(v_i)$. 
			In particular, as $C_{i+1}$ is a face of $D_i$, $F(C_{i+1})$ must also be bounded above by $F(v_i)$. 
			Since $C_{i+1}$ is not in the lower star of $v_i$ (and is in the lower star of $v_{i+1}$), we conclude that that $F({v_{i+1}})$ is strictly less than $F({v_i})$. 
			As a result, all $V$-paths cannot return to the lower star of a vertex once they have left it. 
			
			By construction $V$ has exactly one critical $k$-cell for each vertex $v$ which is a PL critical point of $F$ with index $k$, and a critical $0$-cell $\{*\}$ with $F$-value $-\infty$. 
			Each critical $k$-cell has maximal value given by $F(v)$ for its corresponding vertex $v$. As the vertex $v$ is an index-$k$ PL critical point, $H_k(\CF^{F(v)}, \CF^{F(v)-\epsilon})$ is rank 1  (\cite{GLMas}); this is precisely what is needed for $V$ to be relatively perfect to $F$ (\Cref{def:relperfect}).

		\end{proof}
		
		In summary, the results of this section demonstrate a constructive algorithm for producing a relatively-perfect discrete gradient vector field to $F$ on the cells of $\CF^-$.

	\section{Computational considerations}
	
	\label{sec:computation}
	
	Computational implemention of the algorithm in \Cref{thm:relativelyperfectDGVF} would rely upon identifying whether a vertex in $\CF$ is PL Morse, which involves computing the sign of the gradient on each edge incident to a vertex. 
	Gradients are a local computation, but until now, identifying edges and vertices of this polyhedral complex computationally relied upon an algorithm which globally computes the location and sign sequence of all vertices of $\CF$. 
	This is inefficient if, for example, we wish to follow the local gradient flow and identify reasonable critical cells locally.  In this section, we discuss algorithms which may be used to compute the $\nabla F$-orientation of edges of $\CF$ locally. This computation allows us to construct the pairing from \Cref{thm:relativelyperfectDGVF} locally to a given vertex, including whether this vertex is critical.

	\subsection{Partial Derivatives along Edges}
	
	Here we develop an analytic description of the gradient of $F$ restricted to a cell $C$ of $\CF$, determined by the sign sequence of $C$ (\Cref{lem:cgrad}). This gradient may then be used to explicitly compute a vector in the direction $\overrightarrow{vE}$, that is, from a vertex to an incident edge for any vertex-edge pair in $C$ (\Cref{lem:compute_vE}). By multiplication, we may locally obtain the $\nabla F$-orientation on $E$ (\Cref{cor:compute_gradF}).
	\newpage 
	\begin{lemma}\label{lem:cgrad}
		Let $F$ be a supertransversal neural network given by $$F = G \circ F_m \circ ... \circ F_1$$ Let $C$ be any top-dimensional cell of $\CF$. Then,			
		$$F_{(i)}'|_{C} = \prod_{k=1}^i \relu(\Delta_k (C)) W_k  \quad \textrm{and} \quad \nabla F|_C = W_G \prod_{k=1}^{m} \relu (\Delta_k(C)) W_k$$ 
				
		where $W_k$, $W_G$ are the linear weight matrices of $A_k$ and $G$, respectively, and $\Delta_k(C)$ is the diagonal $n_k\times n_k$ matrix with $s_{kj}(C)$ in each diagonal entry, and $j$ ranges from $1$ to $n_k$. 
	\end{lemma}
	
	\begin{proof}
		By construction $F\big|_{C}$ is affine, and in fact each intermediate $F_{(i)}|_C$ is also affine. Recall that each layer map $F_i$ is given by 
		
		$$F_i = \relu \circ A_i $$
		
		where $A_i$ is an affine function and $\relu$ is the $\max\{0,x\}$ function applied coordinatewise. 
		
		By definition, $s_{ij}(v)=0$ if and only if $\pi_j \circ A_i\circ F_{(i-1)}(v) =0$. 
		This is an affine map. If $A_i \vec{x} = W_i \vec{x} + b_i$ for a weight matrix $W_i$ and bias vector $b_i$, then we note that $F_i\big|_{F_{(i-1)}(C)}$ is given by 
		
		$$F_i\big|_{F_{(i-1)}(C)} (\vec{x})  = \relu( \Delta_i(C)) \left( W_i \vec{x} + \vec{b}_i\right) $$ 
		
		where $\Delta_i(C)$ is the diagonal $n_i\times n_i$ matrix with $s_{ij}(C)$ in each diagonal entry, where $j$ ranges from $0$ to $n_i$ (each of the output entries of $F_i$). 
		This resulting map is affine. By composing these layer maps, we can express $F_{(i)}|_{C}$ as 
		
		$$F_{(i)}|_{C}(\vec{x})=\left(\prod_{k=1}^i \relu(\Delta_k (C)) W_k\right)\vec{x} + \vec{b_i}(C)$$ 
		
		where $\vec{b_i}(C)$ is determined by expanding the composite matrix multiplication (and will ultimately not matter). Resultingly,
		
		$$ F_{(i)}'|_{C} = \left(\prod_{k=1}^i \relu(\Delta_k (C)) W_k\right)$$
		
		Now that we have an equation for $F_{(i)}'|_{C}$. Letting $i=n_m$ (the last layer of $F$) we obtain the total gradient of $F$ given by $F_{m}'$ followed by $W_G$, for $G: \R^{n_m}\to \R$. 
		
		$$\nabla F|_C = W_G \prod_{k=1}^{n_m} \relu(\Delta_k(C)) W_k$$  
		
		where $W_G$ is the weight matrix of the final affine function $G: \mathbb{R}^{n_m} \to \mathbb{R}$, as desired.
		
	\end{proof}\newpage 
	
	We can obtain the gradient of $F$ along edges of $\CF$ by first identifying the direction of $\vec{vE}$, which involves solving an $n_0 \times n_0$ system of equations. 
	
	\begin{lemma} \label{lem:compute_vE}
		Let $F$ be a supertransversal, generic neural network and let $v$ be a vertex of $\CF$. Let $E$ be an edge of $\CF$ incident to $V$. Let $s(v), s(E)$ be sign sequences of $v$ and $E$ respectively. Let $\overrightarrow{vE}$ denote the positive ray of vectors spanned by a vector beginning at $v$ and ending at $x$, a point in $E$. 
		
		Let $\mathcal{I}(v)= \{(i_k, j_k)\}_{k=1}^{n_0}$, consisting of the $(i,j)$ tuples for which $s_{ij}(v) = 0$. Distinguish as $(i_*, j_*)$ the unique $(i,j)$ coordinate for which $s_{i_*, j_*}(E) \neq s_{i_*, j_*}(v)$.  
		
		Denote by $\mathcal{W}(v,C)$ the $(n_0 \times n_0)$ matrix whose $k$th row is given by 
		
		$$ [\mathcal{W}(v,C)]_{k,\cdot} = [ F_{(i_k)}'|_{C}]_{j_k, \cdot } \quad \textrm{where}\quad (i_k, j_k)\in \mathcal{I}(v)$$
		
		Then 
		
		 $$\overrightarrow{vE} = c \cdot s_{i_*j_*}(E) \cdot \mathcal{W}(v,C)^{-1} e_{j_*} $$ 
		 
		 where $c$ is any positive scalar and $e_{j_*}$ is the standard basis vector in $\R^{n_0}$ with a $1$ in the $j_*$ entry and $0$ elsewhere.
		 
	\end{lemma}

	\begin{proof}
		
		We begin by finding a system of equations that $v$ satisfies, determined by its sign sequence.

		Letting $s(v)$ be the sign sequence of $v$, let  $C$ be the cell containing $v$ and $E$ whose sign sequence is obtained by replacing $s(E)$ with $+1$ for all entries that $s(E)=0$. By \Cref{lem:ssprops}, $C$ contains $E$ (and $v$). 
		
		Since $v \in C$ we can express the location of $v$ in $\R^{n_0}$ by the solution of the system of $n_0$ equations given by equations obtained from the $j$th rows of the $i$th intermediate layer maps for each $(i,j)$ pair where $s_{ij}(v)=0$:
		
		\begin{equation}
			[F_{(i)}'|_{C} \vec{v} + \vec{b}_i(C)]_j = 0
			\label{eq:vsystem}
		\end{equation}

		Simultaneously, along every point $x$ in $E$, we have (assuming without loss of generality that $s_{i_*j_*}(E)=1$) that 
		
		$$[F_{(i_*)}'|_{C} x + \vec{b}_{i_*}(C)]_{j_*} > 0$$ 
		
		with all other equations in \Cref{eq:vsystem} satisfied exactly.
		
		Now consider the direction $\vec{x}-\vec{v}$, that is, the direction $\overrightarrow{vE}$. 	
		In each of the $(i,j)$ pairs for which $s_{ij}(v)=0$, we have:

		\begin{align*}
			[ F_{(i)}'|_{C}(\vec{x}-\vec{v})]_j &= 	[F_{(i)}'|_{C} \vec{x} + \vec{b}_{i}(C)]_{j} - [F_{(i)}'|_{C} \vec{v} + \vec{b}_i(C)]_j  \\ 
			& = \begin{cases}
				c & \textrm{with }sign(c)=s_{i_*j_*}(E)\textrm{ when } (i,j)=(i_*, j_*) \\ 
				0 & \textrm{otherwise}
			\end{cases}
		\end{align*}
		
		Selecting any value for $c$ with the appropriate sign gives a $n_0 \times n_0$ system of equations whose solution is a vector in the same direction $\vec{x}-\vec{v}$, that is, $\overrightarrow{vE}$. 
	\end{proof}

	\begin{corollary}\label{cor:compute_gradF}
		Let $F$ be a supertransversal neural network and $v$ be a vertex of $\CF$. Let $E$ be an edge incident to $v$. Let $s(v), s(E)$ be the sign sequences of $v$ and $E$, respectively. Denote by $\partial_{vE}F$ the partial derivative of $F$ in the direction $\overrightarrow{vE}$. 
		
		Then,

		$$\sgn \left(\partial_{vE}F\right) = \sgn\left(\nabla F\big|_{C} \overrightarrow{vE} \right)$$
		
		where $i_*, j_*$ is the layer and neuron, respectively, for which $s_{i*j*}(v)\neq s_{i*j*}(E)$. 
		
	\end{corollary}

	\begin{proof}	
	
		Multiply.
		
	\end{proof}

	\begin{remark}
		Observe that obtaining the sign of $\partial_{vE}F$ involves a total complexity of $n_m+1$ matrix multiplications (followed by ReLU coordinatewise), storing $n_0$ rows of the intermediates to obtain a system of equations to solve for $\overrightarrow{vE}$. This process is comparable in complexity to evaluating $F$ at a point as long as $n_0$ is relatively low. 
	\end{remark}

	\subsection{Following Gradient Flow}

	\label{sec:gradflow}
	
	Here we describe an algorithm which can be used to locally compute a discrete gradient vector field  at a single vertex for a fully-connected, feedforward ReLU neural network $F$ on $\CF$ without having computed all cells of $\CF$. 
	
	In other words, suppose we have a cell $C \in \CF$ with known sign sequence $S(C)$, and we wish to identify a pairing for $C$ in a vector field $V$ compatible with $F$ while relying on a minimal number of evaluations of $F$ or $\partial F_{v\vec{E}}$ at individual points or along individual edges, respectively (as both computations are of similar complexity).
	
	Knowing $S(C)$ and the weights of $F$ gives an explicit set of linear inequalities which bound $C$. 
	The maximal (or minimal) value of $F$ (or any affine function) on $C$ can be identified via linear programming. 
	In practice, applying a solver which uses the simplex algorithm \cite{simplexalg} will quickly identify not only the vertex $v\in C$ where the maximum (or minimum) of $F$ on $C$ is obtained, but also identifies the $n_0$ equations giving its precise location. 
	Setting the signs of $S(C)$ in the entries corresponding to those $n_0$ equations to zero will identify the sign sequence $S(v)$.	
		
	The zero entries of $S(v)$ have an order determined by the existing parameter order of $F$; denote $s_{ij}(v)$ as the sign of $v$ in the $i$th layer and $j$th neuron ($j$th coordinate direction in $\R^{n_i}$). Ordering the entries of $s(v)$ in lexicographic order induces an ordering $\Gamma$ on edges $E$ incident to $v$ by first ordering by which entry of $S(E)$ does not equal that of $S(v)$, and second by whether the entry is negative or positive. 
	
	In order to apply \Cref{lem:local-acyclic-pairing} to $v$ we must evaluate the directional derivative of $F$ along each of the $2n_0$ edges incident to $v$ whose sign sequences can be constructed. This can be done analytically via \Cref{cor:compute_gradF}. If $v$ is regular, then take $e_*$ to be the first edge (with respect to the ordering $\Gamma$) where $\partial F\big|_{e_*} $ is negative but $\partial F\big|_{e_*^{op}}$ is positive; denote its sole additional nonzero entry the $i_*j_*$-th entry.  
	
	Then if $s_{i_*j_*}(C)=0$, we have that $C$ is paired with its coface $D$ obtained by replacing $S(C)$ in index $i_*j_*$ with $s_{i_*j_*}(e_*)$. In the case where $s_{i_*j_*}(C)=s_{i_*j_*}(e_*)$ the cell $C$ is paired with its face $D$ obtained by replacing $s_{i_*j_*}(C)$ with $s_{i_*j_*}(D)=0$. 
	
	If $v$ is critical of index $k$, then take the critical $k$-cell to be the cell obtained by replacing the $k$ zero entries of $s(v)$ that identify the cells in $\st^-(v)$ by $-1$. This  an ordering and labeling of the edges $e_i^-$ for $1 \leq i \leq k$ given by the restriction of $\Gamma$ to these edges. This completes the \textit{Selection Step}, and the algorithm may proceed.
	
	In this way, in order to identify a pairing to which $C$ belongs, we only need to 1) optimize $F$ on $C$, then 2) evaluate $2n_0$ directional derivatives of $F$. The pairings are, in this sense, not determined by the value of $F$ on neighbors of $v$, but instead the choice of ordering of the coordinates in each layer's $\R^{n_i}$ and the relative orientations of their corresponding bent hyperplanes at $v$.

\section{Conclusion}

In this work, we introduce a schematic for translating between a given piecewise linear Morse function on a canonical polyhedral complex and a compatible (``relatively perfect") discrete Morse function. Our approach is constructive, producing an algorithm that can be used to determine if a given vertex in a canonical polyhedral complex corresponds to a piecewise linear Morse critical point, and furthermore an algorithm for constructing a consistent pairing on cells in the canonical polyhedral complex which contain this vertex. However, though we discuss the principles necessary, we leave explicit computational implementation and experimental observations for future work. 

As discussed in \cite{GLMas}, not all ReLU neural networks are piecewise linear Morse, and this is a limitation of our work. Neural networks with ``flat" cells (on which $F$ is constant) are not addressed by our algorithm. This work also defines homological tools which can be used to describe the local change in sublevel set topology at a subcomplex of flat cells, however extensive technical work is needed to provide a direct analog between the cellular topology of $\CF$ and the relevant sublevel set topology. We also leave this to future work. 

We have reason to believe that our proposed algorithm is applicable to any setting in which the star  neighborhoods of the vertices of a PL manifold with the structure of a polyhedral complex are locally combinatorially equivalent to a cross-polytope. The only broad class of functions which we are aware of that satisfies these conditions are ReLU neural networks and similar (for example, leaky ReLU networks or piecewise linear neural networks with activation functions that have several nonlinearities).

We intend that this work be used to develop further theoretical and computational tools for analyzing neural network functions from topological perspectives. 

\label{sec:conclusion}

\section{Acknowledgements}

Many thanks to Eli Grigsby, without whom we may have never put our heads together.

	\bibliographystyle{plain} 
	\bibliography{Discrete_Morse_Canonical_Polyhedra}
\end{document}